\documentclass[12pt,reqno]{amsart}
\usepackage[margin=1in]{geometry} % full-width

\usepackage{times}
\usepackage[english]{babel}
\usepackage[utf8]{inputenc}
\usepackage{csquotes}
\usepackage{cancel}
\usepackage{hyperref}
\usepackage{amsmath,amssymb,amsthm,bbm}
\usepackage{amsfonts}
\usepackage{mathtools}
\usepackage{tabulary}
\usepackage{bm}
\usepackage{bbm}
\usepackage{graphicx, color}
\usepackage{svg}
\usepackage{enumitem}
\usepackage{parskip}
\usepackage{soul}

\usepackage[normalem]{ulem}

\usepackage{collectbox}



\newcommand{\dx}{\,{\rm d} x}

\newcommand{\dy}{\,{\rm d} y}
\newcommand{\dt}{\,{\rm d} t}

\newcommand{\dr}{\,{\rm d} r}
\newcommand{\dxi}{\,{\rm d} \xi}

\newcommand{\dhx}[1]{\,{\rm d} \mathcal{H}^{#1}}
\newcommand{\dH}[2]{\,{\rm d} \mathcal{H}^{#2}(#1)}
\newcommand{\dP}[1]{\,{\rm d} \mathcal{P}^{d-1}_{#1}}
\newcommand{\gammat}{\widetilde{\gamma}}
\newcommand{\gammaSP}{\gamma_{SP} }
\newcommand{\gammaSV}{\gamma_{SV} }
\newcommand{\gammaPV}{\gamma_{PV} }
\newcommand{\gammatPV}{\gammat_{PV}}
\newcommand{\gammatSP}{\gammat_{SP}}
\newcommand{\gammatSV}{\gammat_{SV}}

\newcommand{\T}{\mathbb{T}}
\newcommand{\numberset}{\mathbb}
\newcommand{\R}{\numberset{R}}

\newcommand{\N}{\numberset{N}}

\newcommand{\Z}{\numberset{Z}}

\DeclareMathOperator*{\argmin}{arg\,min}

\newtheorem{theorem}{Theorem}[section]
\newtheorem{definition}{Definition}[section]

\newtheorem{lemma}{Lemma}[section]
\theoremstyle{definition}
\newtheorem{remark}{Remark}[section]

\title[Thresholding energies for anisotropic MCF on inhomogeneous substrate]{Convergence of thresholding energies for anisotropic mean curvature flow on inhomogeneous obstacle}
\author[A. Chiesa] {Andrea Chiesa} 
\address[Andrea Chiesa]{University of Vienna, Faculty of Mathematics  and Vienna School of Mathematics,
                Oskar-Morgenstern-Platz 1, A-1090 Vienna, Austria}
\email{andrea.chiesa@univie.ac.at}
	\urladdr{http://www.mat.univie.ac.at/$\sim$achiesa}

\author[K. Svadlenka] {Karel Svadlenka} 
\address[Karel Svadlenka]{Tokyo Metropolitan University, Graduate School of Science, Department of Mathematics,
                Minami-Osawa 1-1, Hachioji, Tokyo, 192-0397 Japan}
\email{karel@tmu.ac.jp}
	\urladdr{https://karel.fpark.tmu.ac.jp/index\_e.html}


\begin{document}
\begin{abstract}
We extend the analysis by Esedo\={g}lu and Otto (2015) of thresholding energies for the celebrated multiphase Bence-Merriman-Osher algorithm for computing mean curvature flow of interfacial networks, to the case of differing space-dependent anisotropies. In particular, we address the special setting of an obstacle problem, where anisotropic particles move on an inhomogeneous substrate. By suitable modification of the surface energies we construct an approximate energy that uses a single convolution kernel and is monotone with respect to the convolution width. This allows us to prove that the approximate energies $\Gamma$-converge to their sharp interface counterpart.
\end{abstract}

\maketitle

%\tableofcontents

\section{Introduction and Setting}\label{Sec: Intro}

\begin{figure}[ht!]
    \centering
\includegraphics[width=0.65\textwidth]{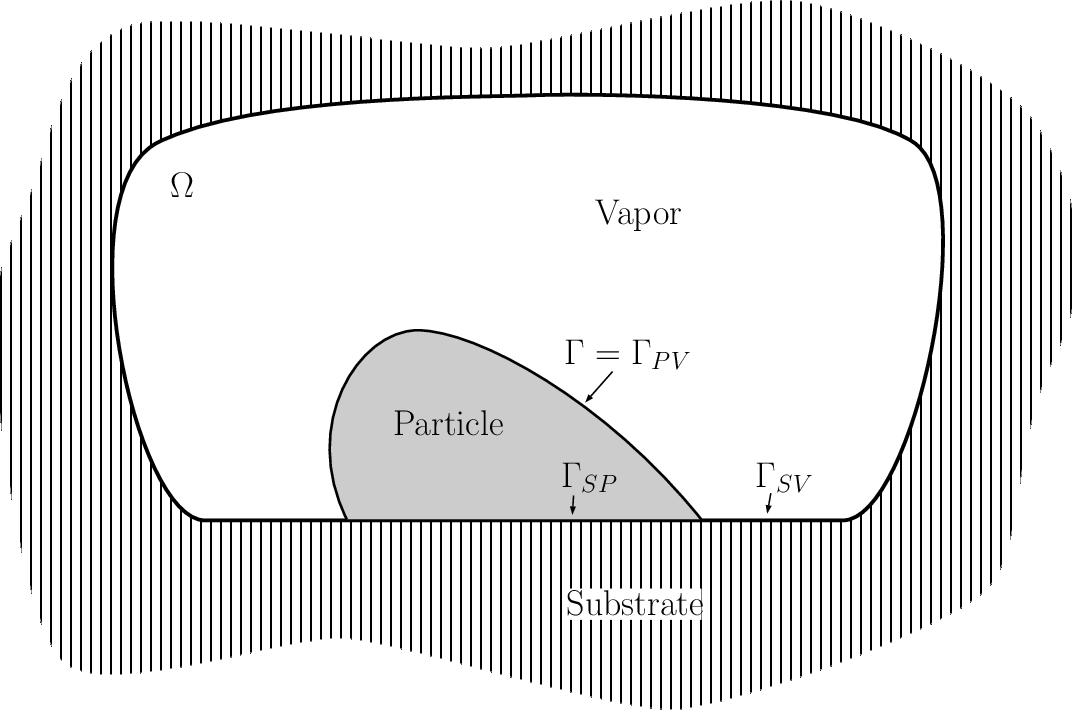}
    \caption{Setup and notation: a particle on a rigid substrate.}
    \label{fig1}
\end{figure}

In this paper, we address the problem of a particle moving on an inhomogeneous substrate according to the $L^2$-gradient flow of the anisotropic interfacial energy.
The setup of the problem is depicted in Figure \ref{fig1}. We let $\Omega$ be a bounded open set of class $C^2$, which is compactly contained in the $d$-dimensional torus $\mathbb{T}^d\coloneqq (\R/\Z)^d$, $d\geq 1$.
We consider the evolution $P(t), \, t\geq 0$, of an open set $P$ representing a particle with fixed mass $m=|P(t)|>0$ in the container $\Omega$. Notice that the assumption of boundedness of $\Omega$ is not restrictive as long as we investigate the time-local behavior of the solution. 
We denote by $V(t) \coloneqq \Omega \setminus \overline{P(t)}$ the complement of $P$ in the container, i.e., the vapor region surrounding the particle, and by $S\coloneqq \mathbb{T}^d\setminus \overline{\Omega}$ the rigid substrate that constrains the motion of $P$. 
Let $\Gamma_{S P}=\partial S \cap \partial P$, $\Gamma_{S V}=\partial S \cap \partial V$, and $\Gamma=\Gamma_{P V}=\partial P \cap \partial V$
be the substrate-particle, substrate-vapor, and particle-vapor interfaces, with surface energy densities or, interchangeably, surface tensions $\gamma_{SP} =\gamma_{SP} (x)$, $\gamma_{SV} =\gamma_{SV} (x)$, and $\psi_{PV}(x,\nu)$, respectively, where $\nu$ is the outer unit normal to $\Gamma$. 
The total energy of the system is given by the anisotropic perimeter
\begin{align*}
    & \int_{\Gamma} \psi_{PV} \dhx{d-1}+\int_{\Gamma_{S P}} \gamma_{SP} \dhx{d-1}+\int_{\Gamma_{S V}} \gamma_{SV} \dhx{d-1} .
\end{align*}
Omitting a constant term depending solely on the data of the problem, we may equivalently consider the energy 
\begin{equation}\label{Physical energy 2}
    \int_{\Gamma}\psi_{PV} \dhx{d-1}+\int_{\partial P\cap \partial S} \sigma\dhx{d-1},
\end{equation}
where $\sigma\coloneqq \gamma_{SP}-\gamma_{SV}$.
Here, $\gamma_{SP}$ and $\gamma_{SV}$ are taken as functions of $x\in \partial S$. We note that this covers the case when the substrate is anisotropic, i.e., possesses normal-dependent energy, since the substrate $S$, and hence also its normal, are fixed in time.

The motion of the particle is given by the $L^2$-gradient flow of the energy \eqref{Physical energy 2}, which amounts to a weighted mean curvature flow interacting with an obstacle.
Namely, the free interface $\Gamma(t)$ evolves in such a way that it fulfills a generalized Herring contact angle condition at the free boundary $\Gamma_{SP} \cap \Gamma_{SV}$, while moving with normal velocity 
\begin{equation}
    \label{grad flow normal velocity}
V = -\mu_{\phi}(\nu)\left(H_\phi-\Lambda_{\phi}\right), 
\end{equation}
where the surface energy density $\phi$, mobility $\mu_{\phi}$, anisotropic mean curvature $H_{\phi}$ and Lagrange multiplier $\Lambda_{\phi}$ are precisely defined in Section \ref{sec_derivation} below.

This kind of interfacial evolution problem appears in several fields of applied science and engineering.
For example, in cell biology it is used as a model to understand cell crawling or patterning \cite{Mohammad2022, Ziebert2016}, while in materials science it assists the development of techniques for coating \cite{Thompson2012} or nanopatterning through the process of solid-state dewetting \cite{Benkouider2015, Bollani2019, Huang2005, Salvalaglio2020}. 
In the mentioned applications, topology changes play a crucial role: cell patterning often involves intercalation of cells, coating is usually achieved by merging of particles, and nanopatterns are often fabricated by fragmentation of films.
To control these processes in terms of scales and interfacial energies, a mathematical approach is required that is able to precisely handle interfacial evolutions across topology changes.

Following the demand in applications, the mathematical and numerical analysis of the moving interface problem is extensive with a rich variety of methods available. It is impossible to provide an exhaustive account of the work done in this field. Therefore, we mention only a few recent results that are closely related to our setting with an obstacle, focusing on the so-called indirect approach that is suitable for dealing with topological singularities.
Bao et al. \cite{Bao2017} analyzed the stationary shapes of anisotropic particles on substrates. 
Mercier and Novaga \cite{Mercier2015} show the existence of a graph solution for the isotropic obstacle problem. Here, the term \emph{isotropic} refers to constant surface energies, independent of the orientation of the interface. 
Kr\"{o}ner et al. \cite{Kroner2021} show the existence of a classical solution to an anisotropic tripod problem. Eto and Giga \cite{Eto2024, Eto2025} follow Bellettini and Kholmatov \cite{Bellettini2018} and analyze a modification of the Almgren-Taylor-Wang approximation scheme for mean-curvature motion of an isotropic droplet with prescribed contact angle, showing its convergence to the viscosity solution of the corresponding level-set equation and proposing a numerical scheme based on the split Bregman method.  
Garcke et al. \cite{Garcke2023} study a phase-field approximation of surface-diffusion motion of anisotropic particles on substrate via anisotropic Cahn-Hilliard equation, showing its relation to a sharp-interface problem based on formal expansions. They also present accurate finite element calculations. 
A similar approach is taken by Hensel and Laux \cite{HenselLaux} who construct BV solutions to isotropic mean curvature flow with constant contact angle via regularization through the Allen-Cahn equation.
On the other hand, Feldman et al. \cite{Feldman2025} address an isotropic problem with linearized surface energy and rate-independent dissipation. 

The above results concern operators different from the anisotropic mean curvature operator that we deal with here.
The same setting as ours has been analyzed only numerically; see, e.g., \cite{Karel, Guo2025}.
Here we present a first step in a rigorous existence result for the anisotropic obstacle problem using the diffuse-interface approach, also called the thresholding method or the BMO algorithm.
This approach is not only powerful analytically, but is also attractive from the viewpoint of numerical solution because it consists of just two fast-computable steps: convolution with an anisotropic kernel and thresholding. 
For the basic problem of mean curvature flow of the boundary of a set $P(t)$ without obstacle and volume constraint and with constant surface energy density equal to $1$, these steps read:
\begin{enumerate}
        \item \emph{convolution}: $\chi^k \coloneqq G_h * \mathbbm{1}_{P^k}$, where $G_h(x)\coloneqq\frac{1}{\sqrt{4 \pi h}^d} e^{-|x|^2/4 h}$ is the scaled Gaussian;
    \item \emph{thresholding}: $P^{k+1} \coloneqq \{x\in \mathbb{T}^d\mid \chi^k(x)> \frac{1}{2}\}$.
\end{enumerate}
They are repeated for $k=0,1,2, \dots$ starting with the characteristic function $\mathbbm{1}_{P^0}$ of the initial set $P^0$.
The interesting point that makes this scheme applicable to a wide range of interface evolution problems is that this thresholding scheme is equivalent to a minimizing movement, i.e., a time-discretization of the gradient flow, for the energy
$$ E_h(u) = \frac{1}{\sqrt{h}} \int (1-u)G_h*u \, \mathrm{d}x , \qquad u\colon \mathbb{T}^d\to [0,1] .$$
It is known that $E_h$ $\Gamma$-converges to the perimeter functional as $h\to 0$ \cite{EsedogluOtto}.

Using the thresholding approach, several problems closely related to ours have been addressed. Xu et al. \cite{Xu2017, Xu2021} investigate isotropic mean curvature flow with obstacle numerically. Esedo\={g}lu, Laux and Otto \cite{EsedogluOtto, LauxLelmi, LauxOtto} develop a theory for the motion of interfacial networks that are isotropic, but each interface is allowed to have a different constant surface tension. Since our obstacle problem can be interpreted as a special type of the three phase problem, this theory can be viewed as a generalization of our problem in the case of isotropic energies.
The authors show that if the surface tensions satisfy a certain triangle inequality condition, the corresponding approximate energies $E_h$ $\Gamma$-converge to the weighted perimeter and the thresholding scheme converges to a weak BV solution.
The thresholding scheme for anisotropic problems without obstacle was first studied in \cite{Ishii1999}, and later \cite{Elsey2017, EsedogluJacobs, EsedogluJacobsZhang} developed a theory of anisotropic kernels $K$ to replace the Gaussian $G$ in the isotropic case.
The analysis focuses on construction of kernels that not only yield the prescribed anisotropy and mobility but also guarantee desirable properties of the resulting thresholding scheme, such as monotonicity and unconditional stability.

When all interfaces in the problem have the same anisotropy, it seems not difficult to extend the results of Laux and Otto \cite{LauxOtto} to the anisotropic network problem. However, it turns out that the analysis of the multiphase problem with differing anisotropies for each interface is challenging. 
The reason is that in such a general situation, it is not possible to show the monotonicity of the approximate energy $E_h$ with respect to $h$, which is then extensively used in the convergence proof. 
Our problem treated in this paper can be viewed as an intermediate step towards understanding the general multiphase problem because we allow for different anisotropies for each of the three types of interfaces but one of them, namely the substrate, is fixed during the evolution.
The presented analysis contributes not only to a new result for the anisotropic mean curvature flow on obstacle but also may provide hints on how to deal with the general multiphase problem.
In particular, in the present paper we construct suitable approximate energies $E_h$ for the anisotropic obstacle problem and prove that their $\Gamma$-limit is the desired interfacial energy \eqref{Physical energy 2}, deferring the proof of convergence of the corresponding thresholding scheme to the following paper.
We achieve this goal under the condition that a strengthened version of the triangle inequality for the surface tensions holds.
This stronger condition comes from the fact that the energy densities are space- and normal-dependent.

Before closing the introduction, we briefly describe the contents of this paper. First, we provide a derivation of the gradient flow equation in the general setting of Finsler geometry, introduced by Bellettini and Paolini in \cite{BellettiniPaolini}.
Here, for a given anisotropy $\phi(x,\nu)$, one considers the anisotropic distance
$$ d_{\phi }(x,y) \coloneqq \inf \left\{ \int_0^1 \phi (c(t),\dot{c}(t))\, \mathrm{d} t \mid c \in W^{1,1}([0,1];\R^d), \, c(0)=x, \, c(1)=y \right\}. $$
It was shown in \cite{BellettiniPaoliniVenturini} that the corresponding anisotropic perimeter of a $C^1$ class set $D$ with normal $\nu_D$ reads $\mathcal{H}^{d-1}_{\phi }(\partial D,\Omega) = \int_{\Omega \cap \partial D} \phi ^o(x,\nu_D) \, \det_{d} \phi \, \dhx{d-1}$, where $\phi^o$ is the dual of $\phi$ and $\det_d\phi$ is the volume element of the anisotropic Hausdorff measure $\mathcal{H}^{d}_{\phi }$.
Our energy in this setting then reads
\begin{equation}\label{def energy finsler}
    E(P) \coloneqq  E_{\phi}(P;\Omega)+E_{\sigma}(P;\overline{S})\coloneqq \int_{\Gamma}\phi^o\operatorname{det}_d\phi\dhx{d-1}+\int_{\partial P\cap \partial S} \sigma\dhx{d-1} .
\end{equation}
In Section \ref{sec_derivation}, under rather weak assumptions on the anisotropic tensions, we calculate the first variation of this energy, revealing the form of the generalized Herring condition that should hold on the free boundary. Subsequently, we derive the corresponding weighted $L^2$-gradient flow and define its solution in the framework of $BV$-functions.

The following Section \ref{sec_thresholding energies} focuses on the particular class of surface tensions with separated variables of the form $\psi_{PV}(x,\nu) = \gamma_{PV}(x)\gamma(\nu)$ and constructs an approximate energy $E_h$ that can be employed in the thresholding scheme. 
As mentioned above, an appropriate nonlocal approximation via convolution with anisotropic kernels is well-known but its straightforward application to our problem would require using different kernels for each type of interface, which would preclude any convergence proof. We deal with this situation by suitably modifying and extending the participating surface tensions to obtain $\gammatPV,\gammatSP,\gammatSV$ and define the approximate energy as 
\begin{equation}\label{energy intro gamma conv}
    E_{h}(u)\coloneqq \frac{1}{\sqrt{h}} \int_{\Omega}\gammatPV  u K_{h} * (\mathbbm{1}_\Omega-u)+\gammatSP u K_{h} * \mathbbm{1}_S +\gammatSV  (\mathbbm{1}_\Omega-u) K_{h} * \mathbbm{1}_S\dx .
\end{equation}
Here $K_h$ is a scaled kernel corresponding to the anisotropy $\gamma(\nu)$, and $u$ is a placeholder for the characteristic function of particle region $P$.

The final Section \ref{bigsec gamma conv} proves that the energy \eqref{energy intro gamma conv} $\Gamma$-converges in the strong $L^1$-topology to \eqref{Physical energy 2} or equivalently to \eqref{def energy finsler} where $ \phi(x,\nu)\coloneqq \sqrt[d-1]{\gammaPV(x)}\,\gamma^o(\nu)$. 
Thanks to the modification of $E_h$ ensuring that all terms are of the form of convolution with the same kernel $K_h$, the proof is allowed to use standard techniques, but there are some difficulties related to the presence of an obstacle and to the fact that we consider $x$-dependent surface tensions.

\section{Derivation of gradient flow equations}
\label{sec_derivation}

In the following sections, we derive the anisotropic mean curvature flow equation \eqref{grad flow normal velocity} as the $L^2$--gradient flow of the energy \eqref{def energy finsler}. We do so in the general Finsler setting introduced in \cite{BellettiniPaolini}, where the surface tension between the particle $P$ and the vapor $V$ is anisotropic and inhomogeneous, i.e., it depends on the normal vector $\nu$ and the spatial variable $x$. On the other hand, the density $\sigma=\sigma(x)$ can be assumed without loss of generality to depend only on $x$ since the substrate $S$ is fixed and, hence, so is the normal vector to $\partial S$ at $x$.

\subsection{Finsler-geometric setting and assumptions}\label{Sec: Admissible anisotropies}

We recall the standard Finsler-geometry setting for anisotropic mean curvature flow introduced in \cite{BellettiniPaolini}.
Let $\Omega\subset \subset \T^d$ be an open, connected set of class $C^2$. 
Let $\phi \colon \overline{\Omega}\times \R^d\to [0,+\infty)$ satisfy
\begin{gather}
\phi  \in C^2(\overline{\Omega}\times (\R^d\setminus \{0\})),\notag\\
\phi ^2(x,\cdot) \text{ is strictly convex for every } x\in \overline{\Omega},\label{condition anisotropy phi}\\
\phi (x, \lambda \nu)=|\lambda| \phi (x, \nu) \text{ for every } x \in \overline{\Omega},\, \nu \in \R^d,\, \lambda \in \R,\notag\\
c_\phi |\nu| \leq \phi (x, \nu) \leq C_{\phi }|\nu| \text{ for every } x \in \overline{\Omega},\,  \nu \in \R^d,\notag
\end{gather}
for some $0<c_{\phi }\leq C_{\phi }$.

Let us denote by $\phi ^o\colon\overline{\Omega}\times \R^d\rightarrow [0,\infty)$ the dual function to $\phi $, i.e., 
\begin{equation*}
    \phi ^o\left(x, \nu^*\right)=\sup \left\{\nu^* \cdot \nu\mid \nu \in B_\phi (x)\right\},
\end{equation*}
where $B_{\phi }(x)=\{\nu\in\R^d\mid \phi (x,\nu)\leq 1\}$. Then $\phi^o$ satisfies the same conditions \eqref{condition anisotropy phi} and $\phi ^{oo}=\phi$ (see \cite{BellettiniPaolini}).

Let us recall some useful results on the anisotropy densities satisfying \eqref{condition anisotropy phi}.

\begin{lemma}[{\cite[Sec. 2.1]{BellettiniPaolini}}]
\label{lemma_phi2}
Let \eqref{condition anisotropy phi} hold. For any $\lambda\neq 0$, $x \in \overline{\Omega}$, and $\nu, \nu^* \in \R^d\setminus\{0\}$ the following identities hold: 
\begin{equation}\label{property anisotropy phi 1}
\phi (x, \nu)=\nabla_{\nu}\phi (x, \nu) \cdot \nu, \quad \phi ^o\left(x, \nu^*\right)=\nabla_{\nu}\phi ^o\left(x, \nu^*\right) \cdot \nu^* ,
\end{equation}
\begin{equation}\label{property anisotropy phi 1.b}
\phi \left(x, \nabla_{\nu}\phi ^o\left(x, \nu^*\right)\right)=\phi ^o\left(x, \nabla_{\nu}\phi (x, \nu)\right)=1 ,
\end{equation}
\begin{equation}\label{property anisotropy phi 1.c}
\phi ^o\left(x, \nu^*\right)\nabla_{\nu}\phi \left(x, \nabla_{\nu}\phi ^o\left(x, \nu^*\right)\right)=\nu^*, \quad \phi \left(x, \nu\right)\nabla_{\nu}\phi ^o\left(x, \nabla_{\nu}\phi \left(x, \nu\right)\right)=\nu.
\end{equation}
In particular, we have, for every $i=1,\dots,d$, 
\begin{equation}\label{property anisotropy phi 2}
    \partial_{x_i}\phi (x,\nabla_{\nu}\phi ^o (x,\nu))+\frac{\nu}{\phi ^o(x,\nu)}\cdot \partial_{x_i}(\nabla_{\nu} \phi ^o (x,\nu))=0,
\end{equation} 
\begin{equation}\label{property anisotropy phi 4}
    \partial_{x_i}\phi (x,\nabla_{\nu}\phi ^o (x,\nu))+\partial_{x_i}\phi ^o\left(x,\frac{\nu}{\phi ^o(x,\nu)}\right)=0.
\end{equation}
\end{lemma}
\begin{proof}
Identities \eqref{property anisotropy phi 1}--\eqref{property anisotropy phi 1.c} can be found in \cite{BellettiniPaolini}.
Formula \eqref{property anisotropy phi 2} follows by differentiating the first equation of \eqref{property anisotropy phi 1.b} and then using the first identity in \eqref{property anisotropy phi 1.c}.
Formula \eqref{property anisotropy phi 4} follows from \eqref{property anisotropy phi 2} combined with the second identity in \eqref{property anisotropy phi 1} differentiated by $x_i$. 
\end{proof}

\subsection{First variation of the energy}

In order to obtain the mean curvature flow equation as the $L^2$-flow of \eqref{def energy finsler}, we compute the first variation of the energy.
In the calculations, we will use the function $\operatorname{det}_d \phi (x) $ and the measure ${\rm d}\mathcal{P}^{n}_{\phi }(x)$, for $n\in \N$, which are defined as
\begin{equation*}
   \operatorname{det}_d \phi (x)\coloneqq \frac{\omega_d}{|B_{\phi }(x)|} \quad \text{ and }\quad {\rm d}\mathcal{P}^{n}_{\phi }(x)
  \coloneqq \operatorname{det}_d \phi (x) \dhx{n}(x) . 
\end{equation*}
Here $|U|$, for a Lebesgue measurable set $U\subset \R^d$, denotes its $d$-dimensional Lebesgue measure, $\omega_d$ stands for the volume of the unit ball in $\R^d$, and $\mathcal{H}^{n}$ for the $n$-dimensional Hausdorff measure. 

\begin{theorem}\label{first variation energy}
    Let $P\subset \Omega$ be an open and bounded subset of $\T^d$ with Lipschitz and piecewise $C^2$ boundary and let $\sigma \colon \partial S \to (0,\infty)$ be Lipschitz continuous. Assume that there exists a function $u\in C^{2}(\T^d)$ such that $P=\{u<0\} \cap \Omega$, $\partial P \cap \Omega =\{u=0\} \cap \Omega$, and $\nabla u\neq 0$ on $\partial^*P$. Let $\xi\in C_c^{\infty }(\T^d;\R^d)$ be an admissible vector field, i.e., such that $\xi|_{\partial S}\cdot \nu_S=0$, where $\nu_S$ denotes the outer unit normal to $S$, and let $\Phi_t\colon\T^d\rightarrow\R^d$ be the vector field solving
    \begin{equation}\label{internal variation}
        \begin{cases}
            \partial_t \Phi_t(x)=\xi(\Phi_t(x))\quad &(t,x)\in (0,\infty)\times \T^d,\\
            \Phi_0(x)=x & x \in \T^d.
            \end{cases}
    \end{equation}
    Then, under the assumptions \eqref{condition anisotropy phi}, letting $P_t\coloneqq \Phi_t(P)$, we have
    \begin{equation*}
        \frac{{\rm d}}{{\rm d}t}E(P_t)\Big|_{t=0}=\int_{\Gamma} H_{\phi }\xi{\cdot} \nu \dP{\phi }(x)+\int_{\partial \Gamma}\phi ^o(x,\nu)\xi_{\tau}{\cdot} b_P\,{\rm d}\mathcal{P}_{\phi }^{d-2}(x)+\int_{\partial (\partial P \cap \partial S)}\!\!\sigma(x)\xi_{\tau}{\cdot} b_S\dhx{d-2}(x),
    \end{equation*}
    where $\Gamma \coloneqq  \partial P \cap \Omega$ is the free surface of $P$, $\nu$ is its outward unit normal, $b_P$ and $b_S$ denote the conormal to the boundary of $\Gamma$ and that of $\partial S\cap \partial P$, respectively, $H_{\phi }$ is the anisotropic mean curvature
    \begin{equation}\label{def H_gamma}
        H_{\phi }(x)=\operatorname{div}(\nabla_{\nu}\phi ^o(x,\nu(x)))+\nabla \log (\operatorname{det}_d \phi  (x)){\cdot} \nabla_{\nu}\phi ^o(x,\nu(x)),
    \end{equation}
    and
    $\xi_{\tau}$ is the tangential component of $\xi$ to $\partial P$, defined as
    \begin{equation}\label{decomposition of xi}
    \xi_{\tau}(x)=\begin{cases}
       \xi(x)- \left(\xi(x)\cdot \frac{\nu}{\phi ^o(x,\nu(x))}\right)\nabla_{\nu}\phi ^o (x,\nu(x))\quad &\text{ if }x\in \Gamma,\\
       \xi(x)- \left(\xi(x)\cdot\nu_S(x)\right)\nu_S(x) \quad &\text{ if } x\in \operatorname{int}(\partial P\cap \partial S).
    \end{cases}
\end{equation}
\end{theorem}
\begin{remark}[Admissibility of $\xi$]
    Notice that thanks to the admissibility of $\xi$, i.e., $\xi|_{\partial S}\cdot \nu_S=0$, and the definition \eqref{internal variation} of $\Phi_t$, it follows that $\Omega$ is a trapping region for the flow $\Phi_t$. Hence,  $P_t\subseteq \Omega$ for all $t\geq 0$, and the perturbation of $P$ is admissible in the sense that it does not penetrate the rigid substrate $S$.
\end{remark}
\begin{remark}[Volume preservation]
\label{remark volume}
    If we want to require volume preservation then we have the additional constraint $\int_{\Gamma}\xi\cdot \nu\dhx{d-1}=0$. Indeed, 
    \begin{equation*}
        0=\left.\frac{{\rm d}}{{\rm d }t}\right|_{t=0}\int_{\Phi_t(P)}\dx =\left.\frac{{\rm d}}{{\rm d }t}\right|_{t=0}\int_{P}|J\Phi_t|\dx =\int_P \operatorname{div}\xi\dx=\int_{\Gamma}\xi\cdot \nu\dhx{d-1},
    \end{equation*}
    where in the last equality we used the divergence theorem and $\xi|_{\partial S}\cdot \nu_S=0$. 
\end{remark}

\begin{proof}[Proof of Theorem \ref{first variation energy}]
    Let us define the mapping $v_t\colon\R^d\rightarrow \R$ by setting $u(x)=v_t(\Phi_t(x))$. Hence, we have
    \begin{equation}
        P_t=\Phi_t(P)=\{v_t<0\}\cap\Omega, \quad \partial P_t \cap\Omega=\{v_t=0\} \cap \Omega, \quad \text{ and } \quad \nabla v_t\neq 0 \text{ on } \partial^* P_t,
    \end{equation}
     for $t$ small enough. The outer normal to $\partial P$ and $\partial P_t$ are then given a.e. by $\nabla u/|\nabla u|$ and $\nabla v_t/|\nabla v_t|$, respectively. Thus, for the anisotropic part of the energy we have
     \begin{equation*}
         E_{\phi }(P_t;\Omega)=\int_{\Gamma}\phi ^o\left(\Phi_t(x),\frac{\nabla v_t}{|\nabla v_t|}(\Phi_t(x))\right)\operatorname{det}_d \phi (\Phi_t(x))\dH{\Phi_t(x)}{d-1}.
     \end{equation*}
In order to compute the $t$-derivative of the above expression, we first recall that 
     \begin{equation*}
         \dhx{d-1}(\Phi_t(x))=\dhx{d-1}(x)+t(\operatorname{div}\xi -\nu \nabla\xi\cdot \nu)\dhx{d-1}(x)+o(t)\dhx{d-1}(x)
     \end{equation*}
    by the change of variable formula on surfaces \cite[Sections 8--9]{Simon}.
Moreover, we claim that 
\begin{equation*}
    \frac{{\rm d}}{{\rm d}t}\left( \Phi_t\frac{\nabla v_t}{|\nabla v_t|}\right)\Big|_{t=0}=-\nu\nabla \xi +(\nu \nabla \xi \cdot \nu)\nu
\end{equation*}
on $\Gamma$. Indeed, differentiating $u(x)=v_t(\Phi_t(x))$ and setting $t=0$, we have
\begin{equation*}
    \frac{{\rm d}}{{\rm d}t}\nabla v_t(\Phi_t(x))|_{t=0}=-\nabla u\nabla \xi,
\end{equation*}
noting that for $t=0$ we have $u(x)=v_{0}(\Phi_0(x))=v_0(x)$.
Thus, we find
\begin{equation*}
    \frac{{\rm d}}{{\rm d}t}\left( \Phi_t\frac{\nabla v_t}{|\nabla v_t|}\right)\Big|_{t=0}=- \frac{1}{|\nabla u|}\nabla u\nabla \xi  +\frac{1}{|\nabla u|^2}\frac{\nabla u\nabla\xi \cdot\nabla u}{|\nabla u|}\nabla u=-\nu \nabla \xi+(\nu  \nabla \xi \cdot\nu)\nu\quad \text{ on }\Gamma,
\end{equation*}
where we used that $\nu=\nabla u/|\nabla u|$.
Combining the equations above we get 
\begin{align}
   \frac{{\rm d}}{{\rm d}t}&E_{\phi }(P_t;\Omega)\Big|_{t=0}={\rm I}+{\rm II}+{\rm III}+{\rm IV}\notag \\
   &\coloneqq \int_{\Gamma}\nabla_x \phi ^o{\cdot} \xi\operatorname{det}_d \phi \dhx{d-1} +\int_{\Gamma} \nabla_{\nu}\phi ^o {\cdot}(-
   \nu\nabla \xi {+}(\nu \nabla \xi{\cdot}  \nu)\nu)\operatorname{det}_d \phi \dhx{d-1} \notag \\
   &\quad+\int_{\Gamma}\phi ^o(\operatorname{div}\xi {-}\nu \nabla\xi {\cdot} \nu)\operatorname{det}_d \phi \dhx{d-1} +\int_{\Gamma}\phi ^o \nabla \operatorname{det}_d \phi {\cdot} \xi \dhx{d-1}
\label{variation four terms}.
\end{align}
Thanks to \eqref{property anisotropy phi 1}, we have
\begin{equation}
\label{variation two three}
    {\rm II}+{\rm III}=\int_{\Gamma} -\nu \nabla \xi\cdot \nabla_{\nu}\phi ^o +\phi ^o\operatorname{div}\xi\dP{\phi }.
\end{equation} 
On the other hand, in order to rewrite ${\rm I}$ in a more convenient way, we decompose $\xi$ as $\xi=\xi_\tau+\xi\cdot\frac{\nu}{\phi ^o}\nabla_\nu\phi ^o$, where $\xi_\tau$ is defined as in  \eqref{decomposition of xi}, and we get 
\begin{align*}
    \left( \frac{\nu}{\phi ^o(x,\nu)}\nabla \xi\right)_i&\stackrel{\eqref{property anisotropy phi 1}}{=}\partial_{x_i} \xi^j_{\tau} \frac{\nu^j}{\phi ^o(x,\nu)} +\partial_{x_i} \left(\xi{\cdot} \frac{\nu}{\phi ^o(x,\nu)}\right)+\left(\xi{\cdot }\frac{\nu}{\phi ^o(x,\nu)}\right)\partial_{x_i}(\partial_{\nu_j}\phi ^o(x,\nu))\frac{\nu^j}{\phi ^o(x,\nu)} \\
    &\stackrel{\eqref{property anisotropy phi 2},\eqref{property anisotropy phi 4}}{=}
    \partial_{x_i} \xi^j_{\tau} \frac{\nu^j}{\phi ^o(x,\nu)} +\partial_{x_i} \left(\xi{\cdot} \frac{\nu}{\phi ^o(x,\nu)}\right)+\left(\xi{\cdot} \frac{\nu}{\phi ^o(x,\nu)}\right)\partial_{x_i}\phi ^o\left(x,\frac{\nu}{\phi ^o(x,\nu)}\right),
\end{align*}
where we sum over repeated indices. Thus, by the positive homogeneity \eqref{condition anisotropy phi} of $\phi ^o$ in the $\nu$-variable, we have shown that
\begin{equation}\label{nabla xi nu over gamma^o}
    \nu\nabla \xi=\nu\nabla \xi_{\tau}+\phi ^o(x,\nu)\nabla\left(\xi{\cdot} \frac{\nu}{\phi ^o(x,\nu)}\right)+\left(\xi{\cdot} \frac{\nu}{\phi ^o(x,\nu)}\right)\nabla_x \phi ^o(x,\nu).
\end{equation}
Hence, we may rewrite ${\rm I}$ as
\begin{align*}
    {\rm I}&\stackrel{\eqref{decomposition of xi}}{=}\int_{\Gamma}\nabla_x \phi ^o (x,\nu){\cdot} \xi_\tau+\left(\xi{\cdot} \frac{\nu}{\phi ^o(x,\nu)}\right)\nabla_x\phi ^o(x,\nu)\cdot \nabla_\nu\phi ^o(x,\nu)\dP{\phi }\\
    &\stackrel{\eqref{nabla xi nu over gamma^o}}{=}
    \int_{\Gamma}\bigg[\nabla_x \phi ^o (x,\nu){\cdot} \xi_\tau+\nu\nabla\xi {\cdot}\nabla_{\nu}\phi ^o(x,\nu)-\nu\nabla\xi_{\tau} {\cdot}\nabla_{\nu}\phi ^o(x,\nu)\\
    &\qquad\quad-\phi ^o(x,\nu)\nabla\left(\xi{\cdot} \frac{\nu}{\phi ^o(x,\nu)}\right){\cdot}\nabla_{\nu}\phi ^o(x,\nu)   \bigg]\dP{\phi }
  \end{align*}
and combining the above terms we find 
\begin{align*}
    {\rm I}+{\rm II}+{\rm III}&\stackrel{\eqref{decomposition of xi}}{=}\int_{\Gamma} \phi ^o(x,\nu)\left(\operatorname{div}\xi_{\tau}+\left(\xi\cdot \frac{\nu}{\phi ^o(x,\nu)}\right)\operatorname{div}(\nabla_{\nu}\phi ^o(x,\nu))  \right)\dP{\phi }\\
    &\quad+\int_{\Gamma}\nabla_x \phi ^o (x,\nu){\cdot} \xi_\tau +\nabla_{\nu}\phi ^o(x,\nu)\nabla \nu\cdot \xi_{\tau}\dP{\phi } \\
    &\quad-\int_{\Gamma}\nabla_{\nu}\phi ^o(x,\nu)\nabla \nu\cdot \xi_{\tau}  +\nu\nabla\xi_{\tau} {\cdot}\nabla_{\nu}\phi ^o(x,\nu) \dP{\phi }.
\end{align*}
For the last term we have
\begin{equation*}
    \int_{\Gamma}  \nabla_{\nu}\phi ^o(x,\nu)\nabla \nu\cdot \xi_{\tau} +\nu\nabla\xi_{\tau} {\cdot}\nabla_{\nu}\phi ^o(x,\nu) \dP{\phi }=
    \int_{\Gamma} \nabla_\nu \phi ^o(x,\nu)\left(\nabla \nu {-}\nabla^\top \nu\right)\cdot \xi_{\tau}\dP{\phi },
\end{equation*}
where we used that $\nu\nabla \xi_\tau+\xi_{\tau}\nabla \nu=0$, which follows by differentiating $\xi_{\tau}\cdot \nu=0$.
Moreover,
\begin{equation*}
    \partial_{x_j} \nu^i= \partial_{x_j} \left(\frac{\partial_{x_i} u}{|\nabla u|}\right)=\frac{\partial^2_{x_i x_j} u}{|\nabla u|}-\frac{\partial_{x_i} u \,\partial_{x_k} u \,\partial^2_{x_j x_k} u}{|\nabla u|^3}=\frac{\partial^2_{x_i x_j} u}{|\nabla u|}-\nu^i\nu^k\frac{\partial^2_{x_j x_k} u}{|\nabla u|} \eqqcolon A_{ij}+B_{ij},
\end{equation*}
where we sum over repeated indices. In view of the symmetry of $A$, we find
\begin{align}
    \nabla_\nu \phi ^o(x,\nu)&\left(\nabla \nu {-}\nabla^\top \nu\right){\cdot} \xi_{\tau}=\nabla_\nu \phi ^o(x,\nu)\left(B{-}B^{\top}\right){\cdot} \xi_{\tau}\notag\\
    &=\partial_{\nu^i}\phi ^o(x,\nu)\left(\nu^j\nu^k\frac{\partial^2_{x_i x_k} u}{|\nabla u|} {-}\nu^i\nu^k\frac{\partial^2_{x_j x_k} u}{|\nabla u|}   \right)\xi_{\tau}^j=-\phi ^o(x,\nu)\nu^k\frac{\partial^2_{x_j x_k} u}{|\nabla u|}  \xi_{\tau}^j\label{antisymmetric term},
\end{align}
where in the last equality we used $\xi_\tau\cdot \nu=0$ and \eqref{property anisotropy phi 1}.
Moreover, notice that 
\begin{equation}\label{second derivative u}
    \partial^2_{x_j x_k} u=\partial_{x_k} (\nu^j|\nabla u|)=|\nabla u|\partial_{x_k}\nu^j+ \nu^j\partial_{x_k}(|\nabla u|).
\end{equation}
Substituting \eqref{second derivative u} into \eqref{antisymmetric term} and again using $\xi_\tau\cdot\nu=0$, we get 
\begin{align*}
    \nabla_\nu \phi ^o(x,\nu)\left(\nabla \nu {-}\nabla^\top \nu\right){\cdot} \xi_{\tau}=-\phi ^o(x,\nu)\nu^k\partial_{x_k} \nu^j  \xi_{\tau}^j=\phi ^o(x,\nu)\nu^k\partial_{x_k}\xi^j_{\tau}\nu^j =\phi ^o(x,\nu)\nu\nabla\xi_{\tau}{\cdot}\nu.
\end{align*}
Hence, we have 
\begin{align*}
    {\rm I}+{\rm II}+{\rm III}
    =&\int_{\Gamma} \phi ^o(x,\nu)\left(\operatorname{div}_{\Gamma}\xi_{\tau}+\left(\xi\cdot \frac{\nu}{\phi ^o(x,\nu)}\right)\operatorname{div}(\nabla_{\nu}\phi ^o(x,\nu))  \right)\dP{\phi }\\
    &+\int_{\Gamma}\nabla_x \phi ^o (x,\nu){\cdot} \xi_\tau +\nabla_{\nu}\phi ^o(x,\nu)\nabla \nu\cdot \xi_{\tau}\dP{\phi } ,
\end{align*}
where $\operatorname{div}_{\Gamma}$ is the surface divergence given by 
$\operatorname{div}_{\Gamma}X\coloneqq \operatorname{div}X-\nu\nabla X\cdot \nu$.
By the divergence theorem on manifolds with boundary, we can rewrite the term involving the divergence of $\xi_{\tau}$ as 
\begin{align*}
    \int_{\Gamma}&\phi ^o(x,\nu)\operatorname{div}_{\Gamma}\xi_{\tau}\dP{\phi }=\int_{\Gamma}\phi ^o(x,\nu)\operatorname{det}_d \phi (x)\operatorname{div}_{\Gamma}\xi_{\tau}\dhx{d-1}\\
    &=-\!\int_{\Gamma}\!\nabla (\phi ^o(x,\nu))\cdot\xi_{\tau}\operatorname{det}_d \phi (x)\dhx{d-1}-\!\int_{\Gamma}\!\phi ^o(x,\nu)\nabla (\operatorname{det}_d \phi (x))\cdot \xi_{\tau}\dhx{d-1} \\
    &\quad+\int_{ \Gamma}\nu(\xi_{\tau} \otimes\nabla(\phi ^o(x,\nu)\operatorname{det}_d \phi (x)))\cdot \nu \dhx{d-1}+\!\int_{\partial \Gamma}\!\phi ^o(x,\nu)\xi_{\tau}{\cdot} b_P\operatorname{det}_d \phi (x)\dhx{d-2}\\
   & = -\!\int_{\Gamma}\!\nabla_x \phi ^o(x,\nu)\cdot \xi_{\tau}+\nabla_{\nu}\phi ^o(x,\nu)\nabla \nu\cdot \xi_{\tau}\dP{\phi }-\int_{\Gamma}\phi ^o(x,\nu)\nabla \log(\operatorname{det}_d \phi (x))\cdot \xi_{\tau}\dP{\phi }\\
   &\quad + \int_{\partial \Gamma}\phi ^o(x,\nu)\xi_{\tau}\cdot b_P \,{\rm d}\mathcal{P}_{\phi }^{d-2},
\end{align*}
where we used that $\nu(\xi_{\tau} \otimes\nabla(\phi ^o(x,\nu)\operatorname{det}_d \phi (x))){\cdot} \nu=\nu^i\xi_\tau^i\partial_{x_j}(\phi ^o(x,\nu)\operatorname{det}_d \phi (x)) \nu^j  =0$ since $\xi_\tau\cdot \nu=0$.
Thus, it follows that 
\begin{align*}
    {\rm I}+{\rm II}+{\rm III}&=\int_{\Gamma} \operatorname{div}(\nabla_{\nu}\phi ^o(x,\nu))\xi\cdot \nu -\phi ^o(x,\nu)\nabla( \operatorname{det}_d \phi (x))\cdot \xi_{\tau}
    \dhx{d-1}\\
    &\quad+ \int_{\partial \Gamma}\phi ^o(x,\nu)\xi_{\tau}\cdot b_P \,{\rm d}\mathcal{P}_{\phi }^{d-2}.
\end{align*}
On the other hand, substituting \eqref{decomposition of xi} in ${\rm IV}$ we have
\begin{align*}
    {\rm IV}=&
\int_{\Gamma}\phi ^o(x,\nu) \nabla (\operatorname{det}_d \phi (x)){\cdot} \xi_{\tau} \dhx{d-1}\\
    &+\int_{\Gamma}\phi ^o(x,\nu)\left(\xi\cdot\frac{\nu}{\phi ^o(x,\nu)}\right) \nabla \log (\operatorname{det}_d \phi (x)){\cdot} \nabla_{\nu}\phi ^o(x,\nu) \dP{\phi }.
\end{align*}
Thus, the combination of the above equalities and the definition \eqref{def H_gamma} of $H_\phi $ yields 
\begin{align*}
    \frac{{\rm d}}{{\rm d}t}E_{\phi }(P_t;\Omega)\Big|_{t=0}&= \int_{\Gamma}H_{\phi }\xi\cdot \nu\dP{\phi }+\int_{\partial \Gamma}\phi ^o(x,\nu)\xi_{\tau}\cdot b_P \,{\rm d}\mathcal{P}_{\phi }^{d-2}.
\end{align*}

If $\partial P\cap \partial S\neq \emptyset$, by using $\xi|_{\partial S}\cdot \nu_S=0$ and the divergence theorem, a simplified version of the above calculation yields 
\begin{align}
    \frac{{\rm d}}{{\rm d}t}E_{\sigma}(P_t;\overline{S})\Big|_{t=0}&=\!\!\int_{\partial S\cap \partial P}\!\!\!\!\!\!\!\nabla \sigma{\cdot} \xi \dhx{d-1} +\!\!\int_{\partial S\cap \partial P}\!\!\!\!\!\!\!\sigma(\operatorname{div}\xi {-}\nu \nabla\xi {\cdot} \nu) \dhx{d-1} \notag \\
   \label{variation sigma}
      &=\!\!\int_{\partial (\partial P \cap \partial S)}\sigma\xi_{\tau}{\cdot} b_S\dhx{d-2}, 
\end{align}
which concludes the proof.
\end{proof}

\subsection{Gradient flow derivation}

In this section we show that the direction of maximal slope of $E$ is a negative multiple of $H_{\phi }\nabla_{\nu}\phi ^o(x,\nu)$. The derivation of the gradient flow also incorporates volume preservation and an arbitrary admissible mobility $\mu_\phi \in C^{0,1}(\mathbb{R}^d;[0,\infty))$ satisfying
\begin{gather}
\label{assumptions mu}
    \mu_{\phi }(\lambda \nu)=|\lambda| \mu_{\phi }(\nu),\quad \text{for every } \, \nu\in \R^d, \; \lambda\in \R,\\
    \mu_{\phi }(\nu)=0 \quad \text{ if and only if } \;\; \nu=0 .\notag
\end{gather}
Let us first introduce the normed space
\begin{equation*}
    L^{2}_{\phi }(\partial P;\R^d)\coloneqq \left\{\xi\in L^2(\partial P; \R^d)\mid  \xi|_{\partial P \cap \partial S }\cdot \nu_S=0\right\}\cap C(\partial P; \R^d),
\end{equation*}
equipped with the norm 
\begin{align*}
\|\xi\|^2_{L^2_{\phi }(\partial P)}\coloneqq \int_{\Gamma}\frac{\phi ^o(x,\nu)^2}{\mu_\phi (\nu)}|\phi (x, \xi(x))|^2\dP{\phi }+
\int_{\partial P\cap \partial S}|\xi(x)|^2\dhx{d-1}.
\end{align*} 
We denote by ${\rm d}E$ the first variation of the energy functional as an element of the dual of the normed space $L^{2}_{\phi }(\partial P;\R^d)$ and by $\langle \cdot,\cdot \rangle$ the duality pairing. 

\begin{theorem}\label{Thm gradient flow}
    Under the assumptions of Theorem \ref{first variation energy} and those in \eqref{assumptions mu}, the problem
    \begin{equation}\label{eq max slope}
        \min\left\{\langle {\rm d}E,\xi\rangle \mid \xi \in L^{2}_{\phi }(\partial P;\R^d), \, \|\xi\|_{L^2_{\phi }(\partial P)}= 1, \int_{\Gamma}\xi{\cdot}\nu\dhx{d-1}=0  \right\}
    \end{equation}
    admits solution provided that
    \begin{equation}\label{Herring angle condition}
\phi ^o(x,\nu)\operatorname{det}_d \phi (x)b_{PS} -(\nabla_\nu\phi ^o(x,\nu)\cdot b_P)\operatorname{det}_d \phi (x)\nu_{PS}+\sigma(x)
b_S=0 \quad  \text{ on } \, \overline{\Gamma} \cap \overline{S},
    \end{equation}
where
\begin{equation*}
    b_{PS}\coloneqq b_P-(b_P{\cdot} \nu_S)\nu_S,\quad \nu_{PS}\coloneqq \nu_{P}-(\nu_P{\cdot} \nu_S)\nu_S,
\end{equation*}
and $\nu_P$ is the normal extended from $\Gamma$ to $\partial \Gamma$ (see Figure \ref{fig:contactpoints}).
    Moreover, if \eqref{Herring angle condition} is satisfied, the solution is given by \begin{equation}\label{thm gradient flow velocity}
    \xi\coloneqq 
        \begin{dcases}
            \frac{\mu_\phi (\nu)}{2\lambda\phi ^o(x,\nu)}\left(H_\phi -\frac{\Lambda}{\operatorname{det}_d\phi }\right)\nabla_{\nu}\phi ^o(x,\nu) \quad &\text{ on } \Gamma,\\
            0 \quad &\text{ on } \partial P\cap \partial S,
        \end{dcases}
    \end{equation}
    where 
    \begin{equation}\label{multiplier volume preservation}
        \Lambda=\frac{\displaystyle\int_{\Gamma}\mu_{\phi }(\nu)H_{\phi }\dhx{d-1}}{\displaystyle\int_{\Gamma}\frac{\mu_{\phi }(\nu)}{\operatorname{det}_d\phi }\dhx{d-1}}
    \end{equation}
    and
    \begin{equation}\label{def lambda small}
        \lambda=-\frac{1}{2}\left(\int_\Gamma \mu_{\phi }(\nu)\left(H_\phi -\frac{\Lambda}{\operatorname{det}_d\phi }\right)^2\dP{\phi }\right)^{1/2}.
    \end{equation} 
\end{theorem} 

This theorem justifies the gradient flow structure of the anisotropic mean curvature flow equation \eqref{grad flow normal velocity}. Indeed, let us consider the velocity $2\lambda \xi$, where $\xi$ is defined in \eqref{thm gradient flow velocity}. On $\Gamma\setminus \partial S$, the tangential component of the velocity does not contribute to the variation of the energy \eqref{def energy finsler} and its normal component is 
\begin{equation}
\label{gradient flow eq}
    V\coloneqq \xi{\cdot }\nu= -\mu_{\phi }(\nu)\left(H_\phi -\frac{\Lambda}{\operatorname{det}_d\phi }\right)
\end{equation}
by \eqref{property anisotropy phi 1}. Thus, we formally obtain \eqref{grad flow normal velocity}.
Notice that in the definition \eqref{thm gradient flow velocity}, $\xi=0$ on $\partial P\cap \partial S$ encodes the fact that the substrate is fixed throughout the evolution. 

\begin{proof}[Proof of Theorem \ref{Thm gradient flow}]
    If $\xi$ is a solution to the minimum problem in the statement of the theorem, then there are $\lambda,\Lambda\in\R$ such that $\xi\in L^{2}_{\phi }(\partial P;\R^d)$ is a stationary point of 
    \begin{equation*}
        F(\xi)\coloneqq \langle {\rm d}E,\xi\rangle-\lambda(\|\xi\|^2_{L^2_{\phi }(\partial P)}-1)-\Lambda\int_{\Gamma} \xi\cdot \nu\dhx{d-1},
    \end{equation*}
    and hence $\frac{{\rm d}}{{\rm d}t}F(\xi+t\zeta)\big|_{t=0}=0$ for every $\zeta\in C^{\infty}_c(\R^d;\R^d)$ with $\zeta|_{\partial S} \cdot \nu_S=0$. By the identification of ${\rm d}E$ thanks to Theorem \ref{first variation energy}, we have 
    \begin{align*}
        0=&\int_{\Gamma} \left(H_\phi  \nu-2\lambda \frac{\phi ^o(x,\nu)^2}{\mu_\phi (\nu)} \phi (x,\xi)\nabla_{\nu}\phi (x,\xi)\right)\cdot \zeta \dP{\phi }-\int_{\partial P\cap \partial S}2\lambda \xi\cdot \zeta \dhx{d-1}\\
        &+\int_{\partial \Gamma} \phi ^o(x,\nu) \zeta_{\tau}\cdot b_P \,{\rm d} \mathcal{P}_{\phi }^{d-2}+\int_{\partial (\partial P \cap \partial S)}\sigma(x) \zeta_{\tau}\cdot b_S\dhx{d-2}-\Lambda\int_\Gamma \zeta\cdot \nu\dhx{d-1}.
    \end{align*}
    Since $\xi|_{\partial S}\cdot\nu_S=0$ and $\zeta$ is arbitrary over $\Gamma$, it follows 
    \begin{gather}
    \label{curvature condition gradient flow1}
        \left(H_\phi -\frac{\Lambda}{\operatorname{det}_d\phi }\right) \nu=2\lambda \frac{\phi ^o(x,\nu)^2}{\mu_{\phi }(\nu)} \phi (x,\xi)\nabla_{\nu}\phi (x,\xi)\quad  \text{ on } \Gamma , \\
    \label{curvature condition gradient flow2}
    \xi=0
\quad  \text{ on } \partial P \cap \partial S .
    \end{gather}
In order to recover the Herring angle condition from the ensuing identity 
\begin{equation}\label{Herring condition gradient flow}
    \int_{\partial \Gamma} \phi ^o(x,\nu) \zeta_{\tau}\cdot b_P \,{\rm d} \mathcal{P}_{\phi }^{d-2}+\int_{\partial (\partial P \cap \partial S)}\sigma(x) \zeta\cdot b_S\dhx{d-2}=0
\end{equation}
for every $\zeta\in C^{\infty}_c(\R^d;\R^d)$ with $\zeta|_{\partial S} \cdot \nu_S=0$, 
let us decompose $\zeta$ in the first integral into its normal and tangential components as in \eqref{decomposition of xi} to obtain
$$ \int_{\partial \Gamma} \phi ^o(x,\nu )\left(\zeta-\frac{1}{\phi ^o(x,\nu)}(\zeta{\cdot} \nu_P)\nabla_\nu\phi ^o(x,\nu)\right)\cdot b_P \,{\rm d} \mathcal{P}_{\phi }^{d-2}+\int_{\partial (\partial P \cap \partial S)} \sigma(x) \zeta{\cdot}b_S\dhx{d-2} = 0 .$$
Since $\zeta$ is arbitrary with $\zeta|_S\cdot\nu_S=0$ , we get \eqref{Herring angle condition}.  

Let us now denote by $T^o_\phi {=}T^o_\phi (x,\cdot)\colon\R^d\rightarrow\R^d$ the map defined by  $T^o_\phi (x,\xi)=\phi ^o(x,\xi)\nabla_{\nu}\phi ^o(x,\xi)$ if $\xi\neq 0$, and  $T^o_\phi (x,\xi)=0$ if $\xi=0$. We notice that 
    \begin{equation*}   T^o_\phi (x,\phi (x,\xi)\nabla_{\nu}\phi (x,\xi))=\phi (x,\xi)\phi ^o(x,\nabla_{\nu}\phi (x,\xi))\nabla_{\nu}\phi ^o(x,\nabla_{\nu}\phi (x,\xi))\stackrel{\eqref{property anisotropy phi 1.b},\eqref{property anisotropy phi 1.c}}{=}\xi ,
    \end{equation*}
    where, in the first equality, we used the 1-homogeneity of $\phi ^o$ and the 0-homogeneity of $\nabla_{\nu}\phi ^o$ in the second variable. Thus, by applying $T^o(x,\cdot)$ to both sides of \eqref{curvature condition gradient flow1} we have 
\begin{equation}\label{expression of xi on Gamma}
        \xi=\frac{\mu_\phi (\nu)}{2\lambda\phi ^o(x,\nu)}\left(H_\phi -\frac{\Lambda}{\operatorname{det}_d\phi }\right)\nabla_{\nu}\phi ^o(x,\nu) \quad \text{ on } \Gamma.
    \end{equation}
    In order to find a solution to \eqref{eq max slope}, we require $\|\xi\|_{L^2_{\phi }(\partial P)}=1$ and $\int_\Gamma\xi\cdot \nu\dhx{d-1}=0$. By \eqref{curvature condition gradient flow1}, \eqref{expression of xi on Gamma}, and by 1-homogeneity of $\phi $, the first condition yields
    \begin{equation*}
        1=\int_{\Gamma}\frac{\phi ^o(x,\nu)^2}{\mu_{\phi }(\nu)}|\phi (x,\xi)|^2\dP{\phi }\stackrel{\eqref{property anisotropy phi 1.b}}{=}\frac{1}{4\lambda^2}\int_{\Gamma}\mu_{\phi }(\nu)\left(H_\phi -\frac{\Lambda}{\operatorname{det}_d\phi }\right)^2\dP{\phi },
    \end{equation*}
    which gives us the value of $\lambda^2$.
    Moreover, by \eqref{expression of xi on Gamma},  \eqref{Herring condition gradient flow}, and \eqref{property anisotropy phi 1}, we have 
    \begin{equation*}
        F(\xi)=\frac{1}{2\lambda}\int_{\Gamma}\left(H_\phi -\frac{\Lambda}{\operatorname{det}_d\phi }\right)^2 \mu_{\phi }(\nu)\dP{\phi }.
    \end{equation*}
    Hence $\xi$ is a solution to \eqref{eq max slope} if the right-hand side is negative, namely if $\lambda$ is given by \eqref{def lambda small}.
    Finally, from
    \begin{equation*}
        0=\int_\Gamma\xi\cdot \nu\dhx{d-1}=\frac{1}{2\lambda}\int_\Gamma \mu_\phi (\nu)\left(H_\phi -\frac{\Lambda}{\operatorname{det}_d\phi }\right)\dhx{d-1},
    \end{equation*}
    we get \eqref{multiplier volume preservation} and conclude the proof.
\end{proof}

\subsection{Notion of solution}
In this section we define the notion of distributional solutions corresponding to the anisotropic mean curvature flow derived in the previous section. There are several notions of weak solutions for the mean curvature flow whose relation is not yet fully understood (see, e.g., the discussion in \cite{FHLS, Laux2024}). Here we introduce a BV-type solution, anticipating the upcoming convergence proof. For that, notice that the energy \eqref{def energy finsler} can be equivalently rewritten for general $u\in BV(\Omega;\{0,1\})$ as
\begin{equation}
\label{Energy general form}
    E(u)=\int_{\Omega} \phi ^o(x,\nu(x))\operatorname{det}_d\phi \,|\nabla u| 
    +\frac{1}{2}\int_{\T^d} \sigma(x) (|\nabla u|{+}|\nabla \mathbbm{1}_{\Omega}|{-}|\nabla(\mathbbm{1}_{\Omega}{-}u)|),
\end{equation}
where 
$\nu(x)$ is the measure-theoretic unit outer normal to $P=\{ x \; | \; u(x) = 1 \}$ defined at points $x$ belonging to the reduced boundary $\partial^*P$.

\begin{definition}[Motion by volume-preserving anisotropic mean curvature flow with obstacle]\label{def notion sol}
Let $T \in(0,\infty)$, $m>0$, and $\Omega$ be an open set, and set $S=\T^d\setminus \overline{\Omega}$. Let 
\begin{equation*}
        u\in BV([0,T]\times \Omega; \{0,1\}) \;\text{ with }\; \sup_{t\in [0,T]} E(u(t))<\infty
\end{equation*}
be trivially extended to zero in $\T^d\setminus\Omega$,
and $\nu\coloneqq -\nabla u/|\nabla u|$ be its measure-theoretic outer normal.
Given a bounded and open  initial configuration $P^0\subset \Omega$ with $E(\mathbbm{1}_{P^0})<\infty$ and $\int_{P_0}\dx =m$, we say that $u$ \emph{moves by (anisotropic) volume preserving curvature flow with obstacle $S$}  
if 
\begin{equation*}
    \int_\Omega u(t,x)\dx =m \qquad \text{ for a.e. } t \in [0,T]
\end{equation*}
 and there exist  $H_\phi ,V\colon [0,T]\times\Omega\rightarrow\R$ with
\begin{equation*}
    \int_{0}^T\int_{\Omega} |H_{\phi }|^2\operatorname{det}_d\phi |\nabla u|\dt<\infty, \quad
    \int_{0}^T\int_{\Omega} |V|^2 \operatorname{det}_d\phi |\nabla u|\dt<\infty ,
\end{equation*}
such that 
\begin{enumerate}
    \item $H_{\phi }$ is the \emph{mean curvature} in the sense that it satisfies
\begin{align}\label{mean curvature volume}
  &\!\!\!\int_0^T\!\!\!\int_{\Omega}\!H_\phi  \xi{\cdot }\nu\operatorname{det}_d\phi |\nabla u|\dt=\int_0^T\!\!\!\int_{\Omega}\!\!\!\left(\phi ^o(x,\nu)\operatorname{div}\xi {-}\nu \nabla\xi {\cdot} \nabla_\nu\phi ^o(x,\nu)\right)\operatorname{det}_d\phi |\nabla u|\dt\notag\\
  &\qquad+\int_0^T\!\!\!\int_{\Omega}\!\!\!\big(\nabla_x \phi ^o(x,\nu){+}\phi ^o(x,\nu) \nabla (\log\operatorname{det}_d\phi (x))\big){\cdot} \xi \operatorname{det}_d\phi |\nabla u|\dt
  \end{align}
for all $\xi \in C^{\infty}\left((0, T) \times\T^d;\R^d\right)$;
\item $V$ is the \emph{normal velocity} in the sense that it satisfies
\begin{align}
\int_{\Omega} \zeta(s,x)  u(s,x)\dx   -\int_{\Omega} \zeta(0,x) \mathbbm{1}_{P^0}(x) \dx  =&\!\!\int_0^s\!\!\int_{\Omega}\!\! \partial_t \zeta(t,x) u(t,x) \dx  \dt\notag\\  
&+\!\!\int_0^s\!\!\int_{\Omega}\!\!  \zeta(t,x) V(t,x)|\nabla u| \dt \label{normal velocity volume}
\end{align}
for all $\zeta \in C^{\infty}\left([0, T] \times\Omega\right)$ and $s\in(0,T]$;
\item $H_{\phi }$ and $V$ satisfy 
\begin{align}
    0=&\!\!\int_0^T\!\!\int_{\Omega}\!\!\left(H_{\phi }-\frac{\Lambda}{\operatorname{det}_d\phi }+\frac{1}{\mu_{\phi }}V\right)\xi \cdot \nu\operatorname{det}_d \phi \,|\nabla u|\dt\notag\\
    &\;+\frac{1}{2}\!\int_0^T\!\!\int_{\T^d}\!\!\left(\sigma(x)(\operatorname{div}\xi {-}\nu \nabla\xi {\cdot} \nu)+\nabla_x \sigma(x){\cdot} \xi\right)\left(|\nabla u|{+}|\nabla\mathbbm{1}_{\Omega}|{-}|\nabla (\mathbbm{1}_{\Omega}-u)|\right)\dt\label{rel V and H volume}
\end{align}
for every $\xi \in C^{\infty}\left((0, T) \times \T^d;\R^d\right)$ such that $\xi|_{(0,T)\times\partial S}\cdot \nu_S=0$, 
where $\Lambda$ is the weighted average of the generalized mean curvature $H_{\phi }$ given by (cf. \eqref{multiplier volume preservation}) 
\begin{equation*}
        \Lambda= \frac{\int_{\Omega}\mu_{\phi }(\nu)H_{\phi }|\nabla u|}{\int_{\Omega}\frac{\mu_{\phi }(\nu)}{\operatorname{det}_d\phi }|\nabla u|}.
    \end{equation*}
\end{enumerate} 
 \end{definition}

\begin{remark}
Papers employing BV-type solutions such as \cite{HenselLaux, LauxStinsonUllrich} additionally require an optimal energy dissipation inequality. Since we are not concerned here with existence of solution, we ignore this condition for simplicity.
\end{remark}

\begin{remark}[No volume preservation]
    In the case of mean curvature flow that does not preserve volume, the correct notion of solution is given by Definition \ref{def notion sol}, setting $\Lambda=0$ in \eqref{rel V and H volume}. 
\end{remark}

\begin{remark}[Spatial independent case without obstacle]
    We remark that this notion of solution corresponds to \cite[Definition 2.8]{LauxStinsonUllrich} when the anisotropies have no dependence on the spatial variable, namely, if $\phi ^o(x,\nu)=\phi ^o(\nu)$, and when, in addition, $\Lambda=0$ and $\Omega=\T^d$. In this setting \eqref{mean curvature volume} and \eqref{rel V and H volume} reduce to
    \begin{align*}
  \int_0^T\int_{\Omega}H_\phi  \xi\cdot \nu|\nabla u|\dt=\int_0^T\int_{\Omega}\left(\phi ^o(\nu)\operatorname{div}\xi {-}\nu \nabla\xi {\cdot} \nabla_\nu\phi ^o(\nu)\right)|\nabla u|\dt,
\end{align*}
\begin{align*}
    \!\!\int_0^T\!\!\int_{\Omega}\!\!\left(H_{\phi }+\frac{1}{\mu_{\phi }}V\right)\xi \cdot \nu|\nabla u|\dt=0,
\end{align*}
respectively.
\end{remark}

\begin{remark}[Spatial independent isotropic case with obstacle]
On the other hand, Definition \ref{def notion sol} corresponds to \cite[Definition 1]{HenselLaux} in the spatial-independent isotropic case, i.e.,   $\phi ^o(x,\nu)=c_0|\nu|$ and $\sigma(x)=\sigma$, and with $\Lambda=0$, $\mu_{\phi}=1$. The combination of \eqref{mean curvature volume} and \eqref{rel V and H volume} then reads as
\begin{align*}
      &c_0\int_0^T\int_{\T^d}\left(\operatorname{div}\xi {-}\nu \nabla\xi {\cdot} \nu\right)|\left(|\nabla u|{+}|\nabla (\mathbbm{1}_{\Omega}{-}u)|{-}|\nabla\mathbbm{1}_{\Omega}|\right)\dt\\
      &\quad+\sigma\!\!\int_0^T\!\!\int_{\T^d}\!\!\left(\operatorname{div}\xi {-}\nu \nabla\xi {\cdot} \nu\right)\left(|\nabla u|{+}|\nabla\mathbbm{1}_{\Omega}|{-}|\nabla (\mathbbm{1}_{\Omega}{-}u)|\right)\dt\\
      &\qquad=-\!\!\int_0^T\!\!\int_{\T^d}\!\! V\xi \cdot \nu\left(|\nabla u|{+}|\nabla (\mathbbm{1}_{\Omega}{-}u)|{-}|\nabla\mathbbm{1}_{\Omega}|\right)\dt.
\end{align*}
\end{remark}

\begin{figure}[ht!]
        \centering
\includegraphics[width=0.6\textwidth]{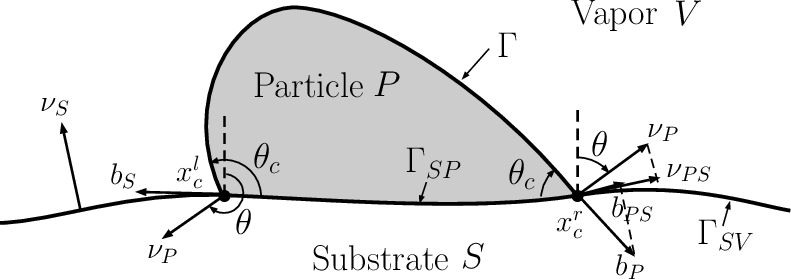}
\caption{Definition of the vectors $\nu_S, \nu_P, \nu_{PS}, b_S, b_P, b_{PS}$ and the contact angles $\theta$ and $\theta_c$.}
\label{fig:contactpoints}
\end{figure}

\begin{remark}[Recovery of classical solution in the smooth case]\label{Remark Recovery of the classical solution}
    In the smooth setting, conditions \eqref{mean curvature volume} and \eqref{rel V and H volume} imply that the normal velocity is proportional to $H_{\phi }$ and that the generalized Herring angle condition holds on the free boundary $\partial\Gamma$. Indeed, doing the same calculations as in the proof of Theorem \ref{first variation energy}, we get
    \begin{align*}
        0=&\int_0^T\int_{\Gamma} \left(\operatorname{div}(\nabla_{\nu}\phi ^o(x,\nu))+\nabla \log (\operatorname{det}_d \phi  (x)){\cdot} \nabla_{\nu}\phi ^o(x,\nu)-\frac{\Lambda}{\operatorname{det}_d \phi }+\frac{1}{\mu_\phi }V\right)\xi{\cdot} \nu \dP{\phi }\dt\\
        &+\int_0^T\int_{\partial \Gamma}\xi_{\tau}{\cdot} \phi ^o(x,\nu)b_P\,{\rm d}\mathcal{P}_{\phi }^{d-2}\dt+\int_0^T\int_{\partial (\partial P \cap \partial S)}\xi_{\tau}{\cdot} \sigma(x)b_S\dhx{d-2}\dt
    \end{align*}
    for every $\xi \in C_c^{\infty}\left((0, T) \times \T^d;\R^d\right)$ such that $\xi|_{(0,T)\times\partial S}\cdot \nu_S=0$. Hence, on one hand, by definition \eqref{def H_gamma} of $H_\phi $ in the smooth setting and the fact that $\xi$ is arbitrary in $(0, T)\times \Omega$ , we find  \eqref{gradient flow eq}. 
    On the other, reasoning as in the proof of Theorem \ref{Thm gradient flow} and using $\xi|_{\partial S} \cdot \nu_S=0$, we recover \eqref{Herring angle condition}. 
    In particular, in the case of flat substrate in $\T^2$, setting $\theta$ to be the angle $\nu_P$ forms with the positive direction of the vertical axis, measured clockwise, we have that 
    \begin{equation*}
        b_{PS} = b_{P}-(b_{P}\cdot\nu_S)\nu_S=\pm (\cos\theta,0), \quad \nu_{PS}=\nu_P-(\nu_{P}\cdot\nu_S)\nu_S=(\sin\theta,0)
    \end{equation*}
    since $\nu_S=(0,1)$, $\nu_P=(\sin \theta,\cos \theta)$ , $b_P=\pm(\cos\theta,-\sin\theta)$, and $b_S=(\pm 1,0)$. Here $+$ or $-$ depends on whether we are considering the right or left contact point, respectively, see Figure \ref{fig:contactpoints}. Hence, assuming $\phi(x,\nu)=\phi^o(\nu)$, \eqref{Herring angle condition} becomes
    \begin{equation*}
    0= \pm \phi(\nu) \cos \theta - (\nabla_\nu\phi(\nu) \cdot b_{P}) \sin \theta \pm \sigma = \pm \widetilde{\phi}\cos\theta \mp \widetilde{\phi}'(\theta)\sin\theta \pm \sigma,
    \end{equation*}
    where $\widetilde{\phi}(\theta)\coloneqq \widetilde{\phi}(\sin\theta,\cos\theta)=\phi(\nu)$. This is the generalized Herring's condition $\widetilde{\phi}\cos\theta - \widetilde{\phi}'(\theta)\sin\theta + \sigma =0$ derived in \cite[Lemma 2.2]{Bao2017}.
    In particular, in the isotropic case, we have $\widetilde{\phi}'(\theta) =0$ and hence $\widetilde{\phi}\cos\theta + \sigma=0$. 
    The contact angle $\theta_c \in [0,\pi)$ measured within the particle $P$ is equal to $\theta$ at the right contact point and to $2\pi-\theta$ at the left contact point (see Figure \ref{fig:contactpoints}). Thus, this identity rewrites $\widetilde{\phi} \cos \theta_c + \sigma = 0$, which is the classical Young's equation.
 \end{remark}

\section{The thresholding scheme and approximate energies}\label{sec_thresholding energies}

In this section, we study an approximation $E_h\colon BV(\Omega;[0,1])\rightarrow [0,\infty)$ of the interfacial energy $E$ defined in \eqref{energy intro gamma conv}
for given $h>0$. To justify the interest in this approximating energy and to illustrate how it arises, let us recall the celebrated thresholding scheme \cite{EsedogluOtto} for the evolution of hypersurfaces by mean curvature flow. For the sake of simplicity, we consider the case of three phases $(\Sigma_i)_{i=1}^3$, disjoint open sets with $\mathbb{T}^d=\cup_{i=1}^3\overline{\Sigma_i}$. Their interfaces $\Gamma_{ij}\coloneqq \overline{\Sigma_i}\cap \overline{\Sigma_j}=\partial \Sigma_i\cap \partial\Sigma_j$, for $i\neq j$, have isotropic and constant surface energy densities $(\gamma_{ij})_{i,j=1}^3 \in \R^{3\times 3}$, where we assume $\gamma_{ii}=0$ and $\gamma_{ij}=\gamma_{ji}$. The energy associated to the partition $(\Sigma_i)_{i=1}^3$ of $\mathbb{T}^d$ is 
\begin{equation*}
    E(\Sigma_1,\Sigma_2,\Sigma_3)\coloneqq \frac{1}{2}\sum_{i,j=1}^3\gamma_{ij}\mathcal{H}^{d-1}(\Gamma_{ij})=\frac{1}{2}\sum_{i,j=1}^3\int_{\Gamma_{ij}}\gamma_{ij}\dhx{d-1},
\end{equation*}
and its $L^2$-gradient flow equations (with mobilities $\mu_{ij}=1/\gamma_{ij}$) read
\begin{equation*}
    V_{ij}=-H_{ij} \quad \text{ on } \Gamma_{ij},
\end{equation*}
where $V_{ij}$ is the normal velocity of $\Gamma_{ij}$ and $H_{ij}$ is its mean curvature.  Additionally, angles at multiple junctions have to satisfy certain conditions. Thanks to the gradient flow structure of these equations, a natural strategy for proving existence of solutions is the minimizing movements scheme: the time interval $[0,T]$ is discretized into $(t_i)_{i=0}^{N_h}$, $t_i=ih$, with time step $h\coloneqq T/N_h$ for some $N_h\in \N$, and an incremental minimum problem of the form
\begin{equation}\label{min probl minimizing movements general}
    u^i\in \argmin_{u\in BV(\mathbb{T}^d;[0,1]^3)}\left\{E_h(u)+\frac{1}{2h}d^2_h(u,u^{i-1})\right\}
\end{equation}
is solved at every step $i=1,\cdots, N_h$, where $u=(u_1,u_2,u_3)$ with $u_1+u_2+u_3=1$ corresponds to $(\mathbbm{1}_{\Sigma_1},\mathbbm{1}_{\Sigma_2},\mathbbm{1}_{\Sigma_3})$. The approximate energy $E_h$ takes the form \cite{EsedogluOtto}
\begin{equation}\label{approximating energy EsedogluOtto}
    E_h(u)\coloneqq \frac{1}{2\sqrt{h}} \sum
_{i,j=1}^3 \int_{\mathbb{T}^d}\gamma_{ij}u_i G_h * u_j\dx,
\end{equation}
where $G(x)= (4\pi)^{-d/2} e^{-|x|^2/4}$ is the Gaussian kernel and the convolution is to be interpreted on $\R^d$. 
The subscript $h$ means rescaling by $\sqrt{h}$: For a function $K\colon \R^d\to\R$ we set
\begin{equation}\label{Kernel scaling}
K_{h}(x)\coloneqq h^{-d / 2} K\left(\frac{x}{\sqrt{h}}\right)\!.
\end{equation}
The induced distance function $d_h$ is defined as $d_h(u,v)\coloneqq E_h(u-v)$. 
Under suitable assumptions including \emph{triangular inequalities}
\begin{equation}\label{triangular inequalities general}
\gamma_{ij}\leq \gamma_{ik}+\gamma_{kj} \quad \text{ for every } i,j,k\in\{1,2,3\},   
\end{equation}
the piecewise constant interpolant in time $\overline{u}_h$, defined as  $\overline{u}_h(t,x)\coloneqq u^i(x)$ for $t\in [t_i, t_{i+1})$, converges in a certain sense to a function $u\in BV(\mathbb{T}^d;\{0,1\}^3)$, solving the gradient flow equation \cite{LauxOtto}. 

One of the main features of the above minimizing movement scheme, making it attractive for numerical solution, is that it can be reinterpreted as a thresholding scheme: For a given initial condition $u^{0}=(\mathbbm{1}_{\Sigma_1^0},\mathbbm{1}_{\Sigma_2^0},\mathbbm{1}_{\Sigma_3^0})$, the solution $\{u^i\}_{i=1}^{N_h}$ to the minimum problems \eqref{min probl minimizing movements general} can be explicitly obtained by repeating the following two steps for $i=1,\dots, N_h$: 
\begin{enumerate}
    \item \emph{convolution}: $\Psi^i_j\coloneqq G_h * \left(\sum_{k=1}^3\gamma_{jk}\mathbbm{1}_{\Sigma^i_k}\right)$, $j\in \{1,2,3\}$; 
    \item \emph{thresholding}: $\Sigma^{i+1}_k\coloneqq \{x\in \mathbb{T}^d\mid \Psi^i_k(x)<\min_{j\neq k} \Psi^i_j(x)\}$.
\end{enumerate}

The definition and analysis of the thresholding scheme for equation \eqref{gradient flow eq} are beyond the scope of this paper and will be addressed in future work. In the following, we consider a first important step in this direction by proposing an appropriate approximation of the energy \eqref{energy intro gamma conv} and showing its $\Gamma$-convergence as $h\rightarrow 0$. 

Motivated by the fact that in the case $\gamma=\gamma(\nu)$ the anisotropy is represented by a suitable kernel $K$ (cf. \eqref{gamma from kernel} below), we consider the case where $\psi_{PV}$ can be multiplicatively split into the product of a space--dependent and a normal--dependent factor  as $\psi_{PV}(x,\nu)=\gammaPV(x)\gamma(\nu)$. To guarantee uniqueness of the decomposition, we assume that $\gamma$ is chosen so that $|B_{\gamma}|=\omega_d$. The relation to the Finsler setting in Section \ref{Sec: Admissible anisotropies} is then given by
$ \phi(x,\nu) = \sqrt[d-1]{\gammaPV(x)} \; \gamma^o(\nu) $. 
The energy \eqref{Physical energy 2} or \eqref{Energy general form} can now, up to a constant, be written as 
\begin{equation}\label{energy gamma conv}
E(u) = \begin{dcases} & \frac{1}{2}\int_{\T^d} \gammaPV (x)\gamma(\nu)\left(|\nabla u|{+}|\nabla(\mathbbm{1}_{\Omega}{-}u)|{-}|\nabla \mathbbm{1}_{\Omega}|\right) \\
& \qquad + \frac{1}{2}\int_{\T^d}\gammaSP (x)\left(|\nabla u|{+}|\nabla \mathbbm{1}_{\Omega}|{-}|\nabla(\mathbbm{1}_{\Omega}{-}u)|\right)\\
&\qquad +\frac{1}{2}\int_{\T^d}\gammaSV (x)\left(|\nabla \mathbbm{1}_{\Omega}|{+}|\nabla(\mathbbm{1}_{\Omega}{-}u)|{-}|\nabla u|\right) \quad \text{for} \; u \in  BV(\Omega;\{0,1\}), \\
& + \infty \qquad\qquad\qquad\qquad \qquad\qquad\qquad\quad \qquad\qquad\text{otherwise}. 
\end{dcases}
\end{equation}

We aim at introducing an associated approximating energy $E_h\colon BV(\Omega;[0,1])\rightarrow [0,\infty)$ in the spirit of \eqref{approximating energy EsedogluOtto} and we notice that in \eqref{approximating energy EsedogluOtto} the kernel used is the Gaussian for all the interfaces. This is an essential ingredient in the proof of a monotonicity formula for the approximating energies (cf. Lemma \ref{lemma gamma conv 2}), which in turn is crucial for the proof of the $\liminf$-inequality part in $\Gamma$-convergence. In the energy \eqref{energy intro gamma conv}, the three interfaces present a different anisotropic behavior, so that the choice of the same kernel $K_h$, say the one associated to $\gamma$, in the approximating energy would result to the wrong energy in the limit. We avoid this issue by modifying the surface tensions $\gammaSP $ and $\gammaSV $ in order to accommodate the presence of anisotropy in the approximate energy, exploiting the fact that the substrate $S$ is fixed. The idea is to define $\gammatSP(x)\coloneqq \gammaSP (x)/\gamma(\nu_S(x))$ and $\gammatSV (x)\coloneqq \gammaSV (x)/\gamma(\nu_S(x))$ for $x\in \partial S$ and to extend them appropriately to $\T^d$.
In Section \ref{Sec: Modified surface tensions}, we provide their construction in detail.
As a result, the energy approximate energy 
\begin{equation}\label{Energy h}
    E_{h}(u)\coloneqq \frac{1}{\sqrt{h}} \int_{\Omega}\gammatPV  u K_{h} * (\mathbbm{1}_\Omega-u)+\gammatSP u K_{h} * \mathbbm{1}_S +\gammatSV  (\mathbbm{1}_\Omega-u) K_{h} * \mathbbm{1}_S\dx 
\end{equation}
$\Gamma$-converges to \eqref{energy gamma conv} in the strong $L^1$-topology, which we prove in Section \ref{bigsec gamma conv}.

Note that the integral in \eqref{Energy h} can be equivalently be taken over $\T^d$ instead of over $\Omega$, if we define $u$ as zero in $\T^d\setminus\Omega$. For technical reasons, in the sequel we will consider the functions $\gammatPV, \gammatSP, \gammatSV$ and $u, \mathbbm{1}_{\Omega}, \mathbbm{1}_S$ appearing in the above definition of $E_h$ to be extended periodically in all coordinate directions of the flat torus $\T^d$, and the convolutions with $K_h$ will be taken over $\R^d$. From the numerical perspective, this does not pose any complications -- in view of the exponential decay of the kernel $K$ (see \eqref{hp K basic} below) we just need to make sure that $\Omega$ is sufficiently far from the boundary of $\T^d$, relative to the magnitude of the smoothing parameter $h$.

\subsection{Assumptions}

Let $\Omega\subset \subset \T^d$ be an open, connected set of class $C^2$. Let $\gamma\colon \R^d\rightarrow [0,\infty)$ be such that 
\begin{gather}
\gamma \in C^2(\R^d\setminus \{0\}), \notag \\
\gamma^2\; \text{ is strictly convex,}
\label{condition anisotropy gamma}
\\
\gamma( \lambda \nu)=|\lambda| \gamma( \nu)\; \text{ for every } \nu \in \R^d,\, \lambda \in \R,
\notag \\
c_\gamma|\nu| \leq \gamma( \nu) \leq C_{\gamma}|\nu|\; \text{ for every } \nu \in \R^d, \notag
\end{gather}
for some  $0<c_{\gamma}\leq C_{\gamma}$.
Moreover, let $\gammaPV\colon \overline{\Omega}\rightarrow (0,\infty)$ and $\gammaSP {,}\gammaSV \colon\partial S\rightarrow (0,\infty)$ be such that 
\begin{gather}
   \gammaPV\in C^{1}(\overline{\Omega})\,  \text{ and }\, \gammaSP ,\gammaSV \in C^{1}(\partial S),
   \label{regularity sigma}
   \\
   c_{s}\leq \gammaSP ,\gammaSV \leq C_s \text{ on }  \partial S,
   \label{boundedness sigma}
   \\
    C_\gamma\gammaPV<\gammaSP +\gammaSV , \quad  \gammaSP <c_\gamma\gammaPV+\gammaSV,\quad 
    \gammaSV < c_\gamma\gammaPV+\gammaSP  \text{ on } \partial S
    \label{condition space gamma triangle ineq}
\end{gather}
for some $0<c_s\leq C_s$.
In view of the decomposition $\psi_{PV}(x,\nu) = \gammaPV(x)\gamma(\nu)$, the condition \eqref{condition space gamma triangle ineq} represents a uniform version of the triangular inequalities \eqref{triangular inequalities general}, and, similarly to \eqref{triangular inequalities general}, is necessary for the lower semicontinuity of the anisotropic energy functional. However, for the extension argument in Section \ref{Sec: Modified surface tensions} below we need to require strict inequalities. \\

\subsection{Modified surface tensions}\label{Sec: Modified surface tensions}

We construct the space-dependent surface tensions $\gammatPV$, $\gammatSP$ and $\gammatSV$, in such a way that they extend the functions $\gammaPV,\gammaSP$ and $\gammaSV$ from their original domains of definition $\overline{\Omega}$ and $\partial S$ to $\T^d$, modify them in a suitable way, so that these modified extensions fulfill
\begin{gather}
\gammatPV ,\gammatSP,\gammatSV \in C^{0,1}( \T^d),
\label{condition space gammat regularity}
\\
\widetilde{c}\leq \gammatPV ,\gammatSP,\gammatSV  \leq \widetilde{C} \quad \text{ in } \T^d  ,
\label{condition space gammat bounds}\\
\gammatPV \leq\gammatSP+\gammatSV , \quad  \gammatSP\leq\gammatPV +\gammatSV , \quad\gammatSV \leq\gammatPV +\gammatSP\quad \text{ in } \T^d,
\label{condition space gammat triangle ineq}
\end{gather}
for some $0<\widetilde{c}\leq \widetilde{C}$.
We remark in particular that the triangular inequalities \eqref{condition space gammat triangle ineq} are meant to hold for every $x\in \T^d$, cf. \eqref{triangular inequalities general} for constant surface tensions. 

To begin with, we define $\gammatPV\colon \T^d\rightarrow (0,\infty)$ as $\gammaPV$ in $\Omega$, and as the solution to the boundary value problem
\begin{equation*}
    \begin{cases}
        \Delta \gammatPV=0 \quad &\text{ in } \T^d\setminus\Omega,\\
        \gammatPV=\gammaPV \quad &\text{ on } \partial\Omega,
        \\
        \gammatPV=\min\limits_{\partial\Omega} \gamma_{PV} \quad &\text{ on }  \partial \T^d,
    \end{cases}
\end{equation*}
in $\T^d\setminus \Omega$. By the regularity theory of the Laplace equation \cite[Thm. 1.8]{Regularity boundary problem} and the maximum principle \cite[Thm. 2.3]{GilbargTrudinger}, $\gammatPV$ satisfies \eqref{condition space gammat regularity} and 
$$ \widetilde{c}_{PV} \leq \gammatPV \leq \widetilde{C}_{PV} \qquad \text{for} \quad \widetilde{c}_{PV} = \min_{\partial\Omega} \gamma_{PV}>0, \; \widetilde{C}_{PV} = \max_{\partial\Omega} \gamma_{PV} .$$

Moreover, since $\Omega\in C^2$ and $\Omega\subset\subset \mathbb{T}^d$, there exists $\delta>0$ such that the open sets $\Omega^+_\delta\subset \Omega$ and $\Omega^-_\delta\subset\mathbb{T}^d\setminus \Omega$, defined as
\begin{equation*}
    \Omega^{\pm}_{\delta}\coloneqq \{x\in \T^d\mid 0<\pm d_{s}(x;\partial\Omega)<\delta\},
\end{equation*}
are of class $C^2$ \cite[Thm. 1]{regularity distance function} and compactly contained in $\mathbb{T}^d$. 
Here we denote by $d_s(\cdot;\partial \Omega)\colon \T^d\rightarrow\R$ the signed distance function from the boundary of $\Omega$, i.e., 
\begin{equation*}
    d_s(x;\partial \Omega)\coloneqq \begin{dcases}
        d(x;\partial \Omega) \quad &\text{ for } x \in\overline{\Omega}, \\
        -d(x;\partial \Omega) \quad &\text{ for } x \in\T^d\setminus \Omega .
    \end{dcases}
\end{equation*}
Let $g^\pm_\delta\colon \Omega^\pm_\delta\rightarrow \R$ be the solutions to the boundary value problem
\begin{equation*}
    \begin{dcases}
        \Delta g^\pm_\delta=0 \quad &\text{ in } \Omega^{\pm}_{\delta},\\
        g^\pm_\delta=\gamma(\nu_{\partial\Omega}) \quad &\text{ on } \partial\Omega,\\
        g^\pm_\delta=c_\gamma \quad &\text{ on } \partial \Omega^{\pm}_{\delta}\setminus \partial \Omega.
    \end{dcases}
\end{equation*}
By \cite[Thm. 1.8]{Regularity boundary problem}, we have
\begin{equation}\label{regularity g_delta}
    g^\pm_\delta\in C^{0,1}(\overline{\Omega^\pm_\delta}).
\end{equation}
Moreover, by the maximum principle \cite[Thm. 2.3]{GilbargTrudinger} and the last condition in \eqref{condition anisotropy gamma}, we have that
\begin{equation}\label{bounds g_delta}
 c_\gamma \leq g^{\pm}_\delta\leq C_\gamma.   
\end{equation}

Let us now extend $\gammaSP$ and $\gammaSV$ to $\T^d$. Let $g_{SP}^\pm, g_{SV}^\pm\colon  \Omega^{\pm}_{\delta} \to (0,\infty)$ be the solutions to 
\begin{equation*}
    \begin{dcases}
        \Delta g^\pm_{SP}=0 \quad &\text{ in } \Omega^{\pm}_{\delta},\\
        g^\pm_{SP}=\gammaSP  \quad &\text{ on } \partial\Omega,\\
        g^\pm_{SP}=\frac{1}{2} C_\gamma\widetilde{C}_{PV} \quad &\text{ on } \partial\Omega^{\pm}_{\delta}\setminus \partial \Omega, 
    \end{dcases}
    \qquad 
    \begin{dcases}
        \Delta g^\pm_{SV}=0 \quad &\text{ in } \Omega^{\pm}_{\delta},\\
        g^\pm_{SV}=\gammaSV  \quad &\text{ on } \partial\Omega,\\
        g^\pm_{SV}=\frac{1}{2} C_\gamma\widetilde{C}_{PV} \quad &\text{ on } \partial\Omega^{\pm}_{\delta}\setminus \partial \Omega, 
    \end{dcases}
\end{equation*}
respectively. By \cite[Thm. 1.8]{Regularity boundary problem} we have
\begin{equation}\label{g_S regularity}
    g^\pm_{SP},g^\pm_{SV}\in C^{0,1}(\overline{\Omega_\delta^\pm}).
\end{equation}
Since the inequalities in \eqref{condition space gamma triangle ineq} are strict, we can assume, up to possibly taking a smaller $\delta$, that 
\begin{equation}\label{triangle ineq g_S}
    C_\gamma\gammatPV<g_{SP}^\pm +g_{SV}^\pm , \quad  g_{SP}^\pm <c_\gamma\gammaPV+g_{SV}^\pm,\quad 
    g_{SV}^\pm < c_\gamma\gammaPV+g_{SP}^\pm \quad \text{ on } \Omega_\delta^\pm.
\end{equation}
We hence define $\gammatSP,\gammatSV \colon \T^d\rightarrow (0,\infty)$ as 
\begin{equation*}
    \gammatSP(x)\coloneqq \begin{dcases}
        \frac{\gammaSP (x)}{\gamma(\nu_{\partial \Omega}(x))} \quad &\text{ if } x \in \partial \Omega,\\
        \frac{g_{SP}^\pm (x)}{g^{\pm}_\delta(x)}\quad &\text{ if } x \in \Omega^{\pm}_\delta,\\
        \frac{C_\gamma \widetilde{C}_{PV}}{2c_\gamma}  \quad &\text{ else}, 
    \end{dcases} 
    \qquad \gammatSV (x)\coloneqq 
    \begin{dcases}
        \frac{\gammaSV (x)}{\gamma(\nu_{\partial \Omega}(x))}\quad &\text{ if } x \in \partial \Omega,\\
        \frac{g_{SV}^\pm (x)}{g^{\pm}_\delta(x)}\quad &\text{ if } x \in \Omega^{\pm}_\delta,\\
        \frac{C_\gamma \widetilde{C}_{PV}}{2c_\gamma} \quad &\text{ else}.
    \end{dcases}
\end{equation*}
By definition, $\gammatPV ,\gammatSP,\gammatSV \in C(\T^d)$ and satisfy the bounds \eqref{condition space gammat bounds} for some $0<\widetilde{c}\leq \widetilde{C}$. Moreover, by \eqref{g_S regularity}, \eqref{regularity g_delta}, and by construction, we have \eqref{condition space gammat regularity}. In order to show the triangular inequalities \eqref{condition space gammat triangle ineq}, let us first consider the case $x\in \overline{\Omega^+_\delta}\cup \overline{\Omega^-_\delta}$. By the triangular inequalities \eqref{triangle ineq g_S} and the bounds  \eqref{bounds g_delta} on $g^\pm_\delta$, it follows 
\begin{equation*}
    \gammatSP(x)+\gammatSV (x)\geq \frac{C_\gamma \gammatPV (x)}{C_\gamma}=\gammatPV (x), \quad |\gammatSP (x)-\gammatSV (x)|\leq \frac{c_\gamma\gammatPV (x)}{c_\gamma}=\gammatPV (x).
\end{equation*}
On the other hand, if $x\in \T^d\setminus (\overline{\Omega^+_\delta}\cup \overline{\Omega^-_\delta})$, then 
\begin{equation*}
    \gammatSP(x)+\gammatSV (x)= \frac{C_\gamma \widetilde{C}_{PV}}{c_\gamma}\geq \gammatPV (x), \quad |\gammatSP (x)-\gammatSV (x)|=0\leq\gammatPV (x).
\end{equation*}

\section{\texorpdfstring{$\Gamma$}{}-convergence of thresholding energies}
\label{bigsec gamma conv}

Here we examine the $\Gamma$-convergence of the approximate energies \eqref{Energy h} toward the original energy functional \eqref{energy gamma conv}. We begin by stating the main result and presenting an outline of its proof. The arguments make use of several lemmas which are then proved in the following subsection.

In this section, we consecutively assume that the surface energy density of the free surface can be written as $\psi_{PV}(x,\nu)=\gammaPV(x)\gamma(\nu)$, where $\gamma(\nu)$ satisfies \eqref{condition anisotropy gamma} and the space-dependent densities $\gammaPV(x),\gammaSP(x),\gammaSV(x)$ satisfy \eqref{regularity sigma}--\eqref{condition space gamma triangle ineq}. Then the modified densities $\gammatPV(x),\gammatSP(x),\gammatSV(x)$ constructed in Section \ref{Sec: Modified surface tensions} meet the conditions \eqref{condition space gammat regularity}--\eqref{condition space gammat triangle ineq}.

Further, let $K\colon \R^d\rightarrow [0,\infty)$ be a kernel associated to the anisotropy $\gamma$, i.e.,
\begin{equation}\label{gamma from kernel}
\gamma(\nu)= \frac{1}{2}\int_0^{+\infty}\!\! r^d  \!\!\int_{\mathbb{S}^{d-1}}\!\! |\xi\cdot\nu| K(r \xi) \dxi \dr ,
\end{equation}
where $\mathbb{S}^{d-1}$ is the $d$-dimensional unit sphere. The kernel $K$ is assumed to satisfy
\begin{gather}
    K(-x)=K(x), \quad K\in L^1(\R^d),\quad K(x)\geq 0\notag,\\ \int_{\R^d}K(x)\dx =1, \quad |x| K(x)\leq c_K K(x/2),\label{hp K basic}
\end{gather}
for some $c_K>0$, so that $|x|K(x)\in L^1(\R^d)$.
Notice that the last condition in \eqref{hp K basic} implies that the kernel decays at least exponentially at infinity.
Moreover, we assume that $K$ is strictly positive in a neighborhood of the origin, namely, there exist positive constants $a,b>0$ such that 
\begin{equation}\label{positive near origin}
    K(x)\geq a \quad \text{ for every } x\in B_b(0).
 \end{equation}
 The condition on the positivity of the kernel $K\geq 0$ can be weakened to sign-changing kernels that possess a positive rearrangement in the sense of \cite[Lemmas 5.6, 5.7]{EsedogluJacobs} since such a rearrangement decreases the energy and does not change the limiting energy.

A question that naturally arises is whether, for a given anisotropy $\gamma$, it is possible to construct a kernel $K$ satisfying the above properties. The problem becomes more involved when one wishes to model also the mobility $\mu_{\phi}$ appearing in Theorem \ref{Thm gradient flow}. We refer to the papers \cite{Elsey2017, EsedogluJacobs, EsedogluJacobsZhang, FuchsLaux, Karel} which address these questions and provide specific methods for the construction of such kernels.

\subsection{Main theorem}\label{Sec: Gamma convercence statement}

The $\Gamma$-convegence result reads as follows.

\begin{theorem}\label{Thm Gamma convergence}
Under the assumptions \eqref{condition anisotropy gamma}--\eqref{condition space gamma triangle ineq} and \eqref{gamma from kernel}--\eqref{positive near origin}, the sequence of functionals $(E_h)_h$, defined as in \eqref{Energy h}, $\Gamma$-converges with respect to the strong $L^1$-topology to the functional $E$, defined in \eqref{energy gamma conv}.
\end{theorem}

The proof of Theorem \ref{Thm Gamma convergence} relies on the following lemmas, whose proofs are presented in Section \ref{Sec: Gamma convergence auxiliary results}. Namely, Lemma \ref{Compactness} provides compactness in the strong $L^1$-topology, Lemma \ref{lemma gamma conv 1} shows the pointwise convergence of $E_h(u)$ to $E(u)$, and Lemma \ref{lemma gamma conv 2} states that the energy $E_h$ is approximately monotonic with respect to $h$.

\begin{lemma}\label{Compactness}
  Consider a sequence of functions $(u_h)_{h\downarrow 0}\subset BV(\Omega,[0,1])$ such that $E_h(u_h)\leq c$ for every $h$, where $E_h$ is as in \eqref{Energy h}, and assume \eqref{condition space gammat regularity}--\eqref{condition space gammat triangle ineq}, \eqref{hp K basic} and \eqref{positive near origin}.
  Then there exists  $u\in BV(\Omega;\{0,1\})$ and a non-relabeled subsequence  $(u_h)_{h\downarrow 0}$ such that $u_h\rightarrow u$ in $L^1(\Omega)$.
\end{lemma}

\begin{lemma}\label{lemma gamma conv 1}
    Let $u\in BV(\Omega;\{0,1\})$ and assume \eqref{condition space gammat regularity}--\eqref{condition space gammat triangle ineq} and \eqref{hp K basic}.
    Then we have 
    \begin{equation}\label{eq: lemma gamma conv 1}
    \lim_{h\rightarrow 0}E_h(u)=E(u).
\end{equation}
\end{lemma}

\begin{lemma}\label{lemma gamma conv 2}
    Assume \eqref{condition space gammat regularity}--\eqref{condition space gammat triangle ineq} and \eqref{hp K basic}.
    Then, for every $u\in BV(\Omega;[0,1])$, every $h>0$ and every $N\in\N_{+}$, we have
    \begin{equation}\label{eq: lemma gamma conv 2}
        E_{N^2 h}(u)\leq (1+cN\sqrt{h}) E_{ h}(u) 
    \end{equation}
    for some $c>0$ independent of $u$, $h$ and $N$.
\end{lemma}

\begin{proof}[Proof of Theorem \ref{Thm Gamma convergence}]$\;$

    $\bullet$ $\limsup$-inequality: We can assume without loss of generality that $u\in BV(\Omega;\{0,1\})$, otherwise $E(u)=\infty$ and the inequality is trivially satisfied. Thanks to Lemma \ref{lemma gamma conv 1}, for $u\in BV(\Omega;\{0,1\})$, we can choose the constant sequence as a recovery sequence and we are done. 

    $\bullet$ $\liminf$-inequality: Let $(u_h)_h\subset BV(\Omega;[0,1])$ with $\liminf_h E_h(u_h)<+\infty$ without loss of generality. By Lemma \ref{Compactness} we can assume $u_h\rightarrow u$ strongly in $L^1(\Omega)$ and $u\in BV(\Omega;\{0,1\})$. 
    Following the strategy of \cite[Proposition 5.5]{EsedogluJacobs}, for $L>0$ let us define
    \begin{equation*}
        K^L(x)\coloneqq K(x)\mathbbm{1}_{B_L}(x),
    \end{equation*}
   and denote by $\gamma^L$ the anisotropy induced by $K^L$ through \eqref{gamma from kernel}. Moreover, we define $E^L$ and $E^L_h$ according to \eqref{energy gamma conv} and \eqref{Energy h}, respectively, with $\gamma$ replaced by $\gamma^L$ and with $K$ replaced by $K^L$. Notice that by positivity of $K$, we have $K^L\leq K$ and $\gamma^L\leq\gamma$, and thus $E_h^L \leq E_h$.
    Fix $h_0\in(0,1)$ and let $N_h\in \N$ be such that $(N_h-1)^2h<h_0\leq N_h^2h$.  Thanks to Lemma \ref{lemma gamma conv 2}, we have 
     \begin{equation*}
         (1+ cN_h \sqrt{h})E_{h}(u_h)\geq E_{N^2_h h}(u_h)\geq E^L_{N^2_h h}(u_h)= E^L_{h_0}(u_h)+R^L_{h,h_0},
    \end{equation*}
    where $R^L_{h,h_0}\coloneqq E^L_{N^2_h h}(u_h)-E^L_{h_0}(u_h)$. Moreover, for fixed $h_0>0$, the functional $E^L_{h_0}$ is continuous with respect to the strong topology of $L^1(\Omega)$ on the space of admissible functions $BV(\Omega;[0,1])$. Hence, we get
    \begin{equation*}
        \lim_{h\rightarrow 0}E_{h_0}^L(u_h)=E_{h_0}^L(u).
    \end{equation*}
    
    Let us show that $\lim_{h\rightarrow 0}|R^L_{h,h_0}|=0$,
    which in turn implies 
    \begin{equation}\label{liminf ineq intermediate}
        (1+c\sqrt{h_0})\liminf_{h\rightarrow 0}E_{h}(u_h)\geq E^L_{h_0}(u)
    \end{equation}
    for every $h_0\in(0,1)$. First, notice
    \begin{align*}
        |R^L_{h,h_0}|&\leq 3c\int_{\Omega}\int_{\R^d}\mathbbm{1}_{B_L}(y)\left|\frac{1}{(N_h^2 h)^{d+1}}K\left(\frac{y}{N_h \sqrt{h}}\right)-\frac{1}{h_0^{d+1}}K\left(\frac{y}{\sqrt{h_0}}\right)\right|\dy\dx \\
        &\leq 3c\int_{\Omega}\int_{B_{L/\sqrt{h_0}}}\frac{1}{h_0}\left|\left(\frac{h_0}{N_h^2 h}\right)^{d+1}K\left(\frac{\sqrt{h_0}}{N_h \sqrt{h}}y\right)-K\left(y\right)\right|\dy\dx .
    \end{align*}
    By density of $C^{\infty}(\R^d)$ in $L^{1}(\R^d)$, for any given $\delta>0$ there exists $\varphi  \in C^{\infty}(\R^d)$ such that 
    \begin{equation*}
       \int_{\R^d}\frac{1}{h_0}\left(\left(\frac{h_0}{N_h^2 h}\right)^{d+1}\left|K\left(\frac{\sqrt{h_0}}{N_h \sqrt{h}}y\right){-}\varphi  \left(\frac{\sqrt{h_0}}{N_h \sqrt{h}}y\right)\right|{+}|\varphi  (y){-}K\left(y\right)|\right)\dy\leq \delta.
    \end{equation*}
    Moreover, since $\Omega$ and $B_{L/\sqrt{h_0}}$ are bounded, and $N^2_h h\rightarrow h_0$ as $h\rightarrow 0$, by uniform continuity we have
    \begin{equation*}
        \lim_{h\rightarrow 0}\int_{\Omega}\int_{B_{L/\sqrt{h_0}}}\frac{1}{h_0}\left|\left(\frac{h_0}{N_h^2 h}\right)^{d+1}\varphi  \left(\frac{\sqrt{h_0}}{N_h \sqrt{h}}y\right)-\varphi  \left(y\right)\right|\dy\dx =0.
    \end{equation*}
    Thus, since $\delta$ is arbitrary, we have $\lim_{h\rightarrow 0}|R_{h,h_0}|=0$. 
    
    Hence, thanks to the pointwise convergence \eqref{eq: lemma gamma conv 1} and \eqref{liminf ineq intermediate}, passing to the limit $h_0\rightarrow 0$, we have 
\begin{equation*}
    \liminf_{h\rightarrow 0}E_{h}(u_h)\geq \lim_{h_0\rightarrow 0}E^L_{h_0}(u)=E^L(u).
\end{equation*}
Letting $L\rightarrow \infty$, we conclude the proof by monotone convergence.
\end{proof}

\subsection{Auxiliary results and proofs}\label{Sec: Gamma convergence auxiliary results}

In this subsection we provide the proofs of Lemmas \ref{Compactness}, \ref{lemma gamma conv 1}, and \ref{lemma gamma conv 2}. First, we state some useful inequalities, which will play a key role in the proofs.\\

\noindent
{\bf Useful inequalities}

\begin{lemma}\label{Lemma useful inequalities}

For any function $v\in L^1(\Omega;[0,1])$ and kernel $K$ satisfying \eqref{hp K basic}, we have 
\begin{gather}
     \!\!\! K_h(y) \!\!\int_{\Omega}|v(x{+}y){-}v(x)| \dx  \dy \leq 2\!\! \int_{\Omega}\!(1{-}v) K_h {*} v \dx +\!\int_{\Omega}    \!\int_{\R^d}\!\!\!\mathbbm{1}_S(x{+}y)K_h(y)v(x)\dy\dx ,\label{ineq 1} \\
     \int_{\Omega}\left|K_h * v-v\right| \dx  \leq \int_{\R^d} K_h(y) \int_{\Omega}|v(x{+}y)-v(x)| \dx  \dy, \label{ineq 2}\\
     \int_{\Omega} v(1-v) \dx  \leq \int_{\Omega}(1-v) K_h * v \dx +\int_{\Omega}\left|K_h * v-v\right| \dx . \label{ineq 3}
\end{gather}
Moreover, if $J$ is a kernel satisfying \eqref{hp K basic} and $ |\nabla J(x)|\leq cJ(x/2)$ for some $c>0$, we have
    \begin{equation}\label{ineq 4}
\int_{\Omega}\left|\nabla\left(J_h * v\right)\right| \dx  \leq  \frac{c}{\sqrt{h}} \int_{\mathbb{R}^d} J_{4h}(y) \int_{\Omega}|v(x{+}y)-v(x)| \dx \dy.
    \end{equation}
\end{lemma}

\begin{proof}$\;$

$\bullet$ To see \eqref{ineq 1}, we first use $K(y)=K(-y)$ to obtain 
    \begin{align*}
        \int_{\Omega}(1{-}v) K_h {* }v \dx &=\int_{\R^d}K_h(y)\int_{\Omega}(1{-}v(x))v(x{+}y) \dx \dy\\
        &=\frac{1}{2}\int_{\Omega}\int_{\R^d}K_h(y)\left((\mathbbm{1}_{\Omega}(x){-}v(x))v(x{+}y){+}(\mathbbm{1}_{\Omega}(x{+}y){-}v(x{+}y))v(x)\right) \dy\dx\\
        &=\frac{1}{2}\int_{\Omega}\int_{\R^d}\mathbbm{1}_\Omega(x{+}y)K_h(y)\left((1{-}v(x))v(x{+}y){+}(1{-}v(x{+}y))v(x) \right)\!\dy\dx.
    \end{align*}
    Hence, \eqref{ineq 1} follows from the elementary inequality 
    \begin{equation*}
        |v'-v|\leq |v'-vv'|+|vv'-v|=v'(1-v)+v(1-v') \quad \text{ for every } v,v'\in [0,1]
    \end{equation*}
    for $v=v(x)$ and $v'=v(x+y)$, and recalling that $K\geq 0$.
Indeed, thanks to the above inequality we get
     \begin{align*}
       & \int_{\Omega}(1{-}v) K_h {*} v \dx  \geq \frac{1}{2} \int_{\Omega}\int_{\R^d}\mathbbm{1}_{\Omega}(x{+}y)K_h(y) |v(x{+}y){-}v(x)| \dy\dx  \\
        & \quad =  \frac{1}{2} \int_{\Omega}\int_{\R^d} K_h(y) |v(x{+}y){-}v(x)| \dy\dx - \frac{1}{2} \int_{\Omega}\int_{\R^d} \mathbbm{1}_{S}(x{+}y)K_h(y) |v(x{+}y){-}v(x)| \dy\dx  \\
        & \quad = \frac{1}{2} \int_{\Omega}\int_{\R^d} K_h(y) |v(x{+}y){-}v(x)| \dy\dx - \frac{1}{2} \int_{\Omega}\int_{\R^d} \mathbbm{1}_{S}(x{+}y)K_h(y) |v(x)| \dy\dx ,
    \end{align*}
which concludes the proof for \eqref{ineq 1}.
    
 $\bullet$ \eqref{ineq 2} follows by Jensen's inequality using the fact that  $K\geq 0$, $\int_{\R^d} K_h=1$ , and $K(-x)=K(x)$. 

$\bullet$ \eqref{ineq 3} follows from the elementary inequality 
    \begin{equation*}
        v(1-v)\leq v'(1-v)+|v'-v| \quad \text{ for every } v,v'\in [0,1]
    \end{equation*}
    applied to $v=v(x)$ and $v'=(K_{h}*v)(x)\in[0,1]$.

    $\bullet$ To prove \eqref{ineq 4}, we first write 
    \begin{equation*}
        \nabla (J_h*v)(x)=\int_{\R^d}\nabla J_h(y)v(x{+}y)\dy=\int_{\R^d}\nabla J_h(y)(v(x{+}y)-v(x))\dy,
    \end{equation*}
    since by $J(x)=J(-x)$ it follows that $\int_{\R^d}\nabla J_h=0$. Hence, we have 
    \begin{equation*}
        \int_{\Omega}|\nabla ( J_h*v)|\dx  \leq \int_{\R^d} |\nabla J_h(y)|\int_{\Omega}|v(x{+}y)-v(x)|\dx \dy.
    \end{equation*}
  Recalling $|\nabla J(x)|\leq c J(x/2)$, we have
  \begin{equation*}
      |\nabla J_h(y)|=h^{-\frac{d{+}1}{2}}|\nabla J(y/\sqrt{h})| \leq ch^{-\frac{d{+}1}{2}} J(y/2\sqrt{h} \leq c h^{-\frac{1}{2}}J_{4h}(y), 
  \end{equation*}
  which in combination with the previous inequality yields \eqref{ineq 4}.
\end{proof}

\noindent
{\bf Proofs of auxiliary lemmas}

\noindent
We first prove the approximate monotonicity of the energy $E_h$ defined in \eqref{Energy h}.

\begin{proof}[Proof of Lemma \ref{lemma gamma conv 2}]
    Fix $u\in  BV(\Omega;[0,1])$ and $h>0$. 
In order to get neat formulas, in the following we will use the new symbols $u^1\coloneqq u$, $u^2\coloneqq \mathbbm{1}_\Omega-u$, and $u^3\coloneqq \mathbbm{1}_S$. Moreover, we put $\gamma_{12}\coloneqq\gamma_{21}\coloneqq\gammatPV $, $\gamma_{13}\coloneqq\gamma_{31}\coloneqq\gammatSP$, 
$\gamma_{23}\coloneqq\gamma_{32}\coloneqq\gammatSV $, and $\gamma_{11}\coloneqq\gamma_{22}\coloneqq\gamma_{33}=0$, and recall that all these functions are periodically extended to $\R^d$. The energy \eqref{Energy h} can be then equivalently rewritten according to this notation as
\begin{equation*}
    E_h(u)=\frac{1}{2\sqrt{h}}\sum_{i,j=1}^3\int_{\T^d}\gamma_{ij}u^iK_h*u^j,
\end{equation*}
since  $K$ is symmetric.
As will be seen below, the inequality \eqref{eq: lemma gamma conv 2} follows by reiterating the formula
\begin{align}
    {\sqrt{h''}}E_{h''}(u)\leq &\sqrt{h}E_{h}(u)+\sqrt{h'}E_{h'}(u)\notag\\
    &+\sum_{i,j=1}^3\int_{\R^d} K(y)\int_{\T^d}u^i(x)u^j(x{-}\sqrt{h}y)|\gamma_{ij}(x)-\gamma_{ij}(x{+}\sqrt{h'}y)|\dx\dy ,\label{eq: lemma gamma conv 2 multiphase K reduced}
\end{align}
where $\sqrt{h''}=\sqrt{h}+\sqrt{h'}$. 
    For notational simplicity, let us define, for fixed $y\in\R^d$,
    \begin{equation*}
       F(\sqrt{h})\coloneqq\sum_{i, j=1}^3 \int_{\T^d}\gamma_{ij}(x)u^i(x)u^j(x-\sqrt{h}y)\dx , 
    \end{equation*}  
     so that $2\sqrt{h}E_h(u)=\int_{\R^d} K F(\sqrt{h})\, \dy$. By change of variable $\tilde{x}=x{+}\sqrt{h'}y$ in the last term of the right-hand side of \eqref{eq: lemma gamma conv 2 multiphase K reduced}, to show \eqref{eq: lemma gamma conv 2 multiphase K reduced} it suffices to prove that for every $y\in \R^d$ 
    \begin{equation}\label{eq: lemma gamma conv 2 multiphase K reduced 1}
        F(\sqrt{h''})-F(\sqrt{h})-F(\sqrt{h'})\leq \sum_{i,j=1}^3\int_{\T^d}\!\!u^i(x')u^j(x'')|\gamma_{ij}(x')-\gamma_{ij}(x)| \dx ,
    \end{equation}
    where $x'\coloneqq x-\sqrt{h'}y, x''\coloneqq x-\sqrt{h''}y$. 
    First, notice that in view of $u^1+u^2+u^3\equiv 1$ we can write for $i\neq j$
    \begin{align*}
        &u^i(x)u^j(x'')-u^i(x)u^j(x')-u^i(x')u^j(x'')\\
        &\qquad=\sum_{k=1}^3 \left(u^i(x) u^k(x')u^j(x'')-u^i(x)u^j(x') u^k(x'')- u^k(x)u^i(x')u^j(x'')\right).
    \end{align*}
Regarding the terms where $k\in\{i,j\}$ in the right-hand side, we have 
\begin{align*}
&u^i(x) u^i(x')u^j(x'')-u^i(x)u^j(x') u^i(x'')- u^i(x)u^i(x')u^j(x'')\\
&\qquad+u^i(x) u^j(x')u^j(x'')-u^i(x)u^j(x') u^j(x'')- u^j(x)u^i(x')u^j(x'')\\
&\quad=-u^i(x)u^j(x') u^i(x'')-u^j(x)u^i(x') u^j(x'')\leq 0.
    \end{align*}
Hence, we have found that 
\begin{align*}
        &u^i(x)u^j(x'')-u^i(x)u^j(x')-u^i(x')u^j(x'')\\
        &\qquad\leq \sum_{k\neq i,j} \left(u^i(x) u^k(x')u^j(x'')-u^i(x)u^j(x') u^k(x'')- u^k(x)u^i(x')u^j(x'')\right).
    \end{align*}
Multiplying both sides by $\gamma_{ij}(x)$, summing over $i,j$, integrating over $\T^d$, and adding and subtracting the term $\sum_{i,j}\int_{\T^{d}}\gamma_{ij}(x')u^i(x')u^j(x'')\dx$, we obtain
\begin{align*}
    &F(\sqrt{h''})-F(\sqrt{h})-F(\sqrt{h'})\\
    &\qquad\leq \sum_{\substack{i,j=1\\ i\neq j}}^3\sum_{k\neq i,j}\int_{\T^d}\gamma_{ij}(x)\left(u^i(x) u^k(x')u^j(x'')-u^i(x)u^j(x') u^k(x'')- u^k(x)u^i(x')u^j(x'')\right)\dx\\
    &\qquad \qquad +\sum_{\substack{i,j=1\\ i\neq j}}^3\int_{\T^d}\!\!u^i(x')u^j(x'')|\gamma_{ij}(x')-\gamma_{ij}(x)|\dx.
\end{align*}
Using the triangular inequalities \eqref{condition space gammat triangle ineq}, we have
\begin{align*}
    &\sum_{\substack{i,j=1\\ i\neq j}}^3\sum_{\substack{k=1\\ k\neq i,j}}^3\int_{\T^d}\gamma_{ij}(x)\left(u^i(x) u^k(x')u^j(x'')-u^i(x)u^j(x') u^k(x'')- u^k(x)u^i(x')u^j(x'')\right)\dx\\
    &\qquad \leq \sum_{\substack{i,j=1\\ i\neq j}}^3\sum_{\substack{k=1\\ k\neq i,j}}^3\int_{\T^d}\gamma_{ik}(x)u^i(x) u^k(x')u^j(x'')+\gamma_{kj}(x)u^i(x) u^k(x')u^j(x'')\dx\\
    &\qquad \qquad -\sum_{\substack{i,j=1\\ i\neq j}}^3\sum_{\substack{k=1\\ k\neq i,j}}^3\int_{\T^d}\gamma_{ij}(x)u^i(x)u^j(x') u^k(x'')+ \gamma_{ij}(x)u^k(x)u^i(x')u^j(x'')\dx\leq 0.
\end{align*}
Consequently, we have shown
\eqref{eq: lemma gamma conv 2 multiphase K reduced 1} and, in turn, \eqref{eq: lemma gamma conv 2 multiphase K reduced}.

Reiterating \eqref{eq: lemma gamma conv 2 multiphase K reduced}, we get, for $N\in \N_+$, 
\begin{equation}\label{gamma conv 2: inequality with error term}
    E_{N^2h}(u)\leq E_h(u)+ C_{N,h}^\gamma(u),
\end{equation}
where 
\begin{equation*}
   C_{N,h}^\gamma(u)\coloneqq \frac{1}{N\sqrt{h}}\sum_{i,j=1}^3\int_{\R^d} K(y)\int_{\T^d}u^i(x)u^j(x{-}\sqrt{h}y)\sum_{n=1}^{N-1}|\gamma_{ij}(x){-}\gamma_{ij}(x{+}n\sqrt{h}y)|\dx\dy. 
\end{equation*}
Thus, it remains to bound the last term of the right-hand side.
Using the lipschitzianity \eqref{condition space gammat regularity} of $\gamma_{ij}$ and the bound $|x|K(x)\leq c_KK(x/2)$, we find 
    \begin{align*}
        C_{N,h}^\gamma(u) &\stackrel{\eqref{condition space gammat regularity}}{\leq} c\frac{N \sqrt{h}}{\sqrt{h}}\sum_{i,j=1}^3\int_{\R^d}|y| K(y)\int_{\T^d}u^i(x)u^j(x-\sqrt{h}y)\dx\dy\\
        &\stackrel{\eqref{hp K basic}}{\leq}  c\frac{N \sqrt{h}}{\sqrt{h}}\sum_{i,j=1}^3\int_{\R^d}
        K\left(\frac{y}{2}\right)\int_{\T^d}u^i(x)u^j(x-\sqrt{h}y)\dx\dy\stackrel{\eqref{condition space gammat bounds}}{\leq}cN\sqrt{h} E_{4 h}(u) .
    \end{align*}
Applying \eqref{gamma conv 2: inequality with error term} with $N=2$, we have 
    \begin{align*}  
    C_{N,h}^\gamma(u)&\stackrel{\eqref{gamma conv 2: inequality with error term}}{\leq} cN \sqrt{h}\left(E_h(u){+}\frac{1}{2\sqrt{h}}\sum_{i,j=1}^3\!\int_{\R^d} \!\!\!\!K(y)\!\!\!\int_{\T^d}\!\!\!\!\!u^i(x)u^j(x{-}\sqrt{h}y)|\gamma_{ij}(x){-}\gamma_{ij}(x{-}\sqrt{h}y)|\right)\dx\dy\\
    &\stackrel{\eqref{condition space gammat bounds}}{\leq} 
        cN \sqrt{h}\left(E_h(u){+}\frac{c}{\sqrt{h}}\sum_{i,j=1}^3\!\int_{\R^d} K(y)\int_{\T^d}u^i(x)u^j(x{-}\sqrt{h}y)\right)\dx\dy\stackrel{\eqref{condition space gammat bounds}}{\leq} c N\sqrt{h} E_h(u),
    \end{align*}
thus concluding the proof.
\end{proof}

We now prove the compactness result.
\begin{proof}[Proof of Lemma \ref{Compactness}]
     Let us define the kernel $J$ as 
    \begin{equation*}
        J(x)\coloneqq \begin{cases}
            ac_J\left(1-\frac{|x|}{b}\right) &\text{ if } x\in B_b(0),\\
            0 &\text{ otherwise,}
        \end{cases}
    \end{equation*}
    where $a,b$ are defined in \eqref{positive near origin} and $c_J$ is such that $\int J = 1$. Notice that $J$
    satisfies \eqref{hp K basic} and 
    $ |\nabla J(x)|\leq 2J(x/2)$. Moreover, we have that  $\frac{1}{c_J}J(x)\leq K(x)$ and $\widetilde{c}\leq \gammatPV (x)$ for every $x\in \R^d$ thanks to \eqref{positive near origin} and \eqref{condition space gammat bounds}, and thus 
    \begin{equation}
        0\leq  \frac{\widetilde{c}}{c_J\sqrt{h}}\int_{\Omega}u_h J_h*(\mathbbm{1}_{\Omega}-u_h)\dx  \leq \frac{1}{\sqrt{h}}\int_{\Omega}\gammatPV u_h K_h*(\mathbbm{1}_{\Omega}-u_h)\dx \leq E_h(u_h)\leq c \label{eq compactness distance characteristic function preliminary 1}.
    \end{equation}
    Hence, since $ |\nabla J(x)|\leq cJ(x/2)$, we can apply \eqref{ineq 4} and we get
    \begin{align}
        \int_{\Omega}|\nabla (J_h*u_h)|\dx  &\stackrel{\eqref{ineq 4}}{\leq} \frac{c}{\sqrt{h}} \int_{\mathbb{R}^d} J_{4h}(y) \int_{\Omega}|u_h(x+y)-u_h(x)| \dx \dy \notag\\
        &\leq \frac{c}{\sqrt{h}} \int_{\mathbb{R}^d} K_{4h}(y) \int_{\Omega}|u_h(x+y)-u_h(x)| \dx \dy \notag\\
        &\stackrel{\eqref{ineq 1}}{\leq} \frac{ c}{\sqrt{h}}\left( \int_{\Omega}(1{-}u_h) K_{4h} {*} u_h \dx  \,{+} \!\int_{\Omega}\int_{\R^d}\! \mathbbm{1}_{S}(x{+}y)K_{4h}(y) u_h(x) \dy\dx\right)\! \label{ineq compactness 1}.
    \end{align}
   Notice that 
    \begin{align}
  \frac{1}{\sqrt{h}}\!\int_{\Omega}\!\int_{\R^d} \!\!\mathbbm{1}_{S}(x{+}y) K_{4h}(y) &u_h(x) \dy\dx = \frac{1}{\sqrt{h}}\!\int_{\Omega}\!\int_{\R^d}\!\!u_h(x)\frac{1}{(4h)^{d/2}}K\left(\frac{y}{2\sqrt{h}}\right)\mathbbm{1}_S(x+y)\dy\dx \notag\\
       & \leq \frac{1}{\sqrt{h}}\!\int_{\Omega}\!\int_{\R^d}\!\mathbbm{1}_{\Omega}(x)\frac{1}{(4h)^{d/2}}K\left(\frac{y}{2\sqrt{ h}}\right)\mathbbm{1}_S(x+y)\dy\dx  \notag \\
       &= \!\int_{\Omega}\!\int_{\R^d}\!\!\mathbbm{1}_{\Omega}(x)K(y)\frac{\mathbbm{1}_S(x{+}2\sqrt{h}y){-}\mathbbm{1}_S(x)}{\sqrt{h}}\dy\dx ,\label{ineq term R^d minus A_x }
   \end{align}
    which is bounded since $\Omega$ has finite perimeter and $|x|K(x)\in L^1$.
   Indeed, rewriting in polar coordinates, we have 
   \begin{align*}
\Bigg|\!\int_{\Omega}\!\int_{\R^d} &\mathbbm{1}_{S}(x{+}y)K_{4h}(y) \frac{u_h(x)}{\sqrt{h}} \dy\dx \Bigg|\\
  &\leq
  \!\int_0^\infty\!\! \!\!\int_{\mathbb{S}^{d{-}1}}\!\!\int_{\T^d}\!
  \!r^d K\left(r\xi\right)\left|\frac{\mathbbm{1}_S(x{+}2r\sqrt{h}\xi){-}\mathbbm{1}_S(x)}{r\sqrt{h}}\right|\dx \dH{\xi}{d-1}\dr\\
  &\leq
  \!\int_0^\infty\!\! \!\!\int_{\mathbb{S}^{d{-}1}}\!\!\int_{\T^d}\!
  \!r^d K\left(r\xi\right)\left|\frac{\mathbbm{1}_\Omega(x{+}2r\sqrt{h}\xi){-}\mathbbm{1}_\Omega(x)}{r\sqrt{h}}\right|\dx \dH{\xi}{d-1}\dr\\
  &\leq \!\int_0^\infty\!\! \!\int_{\mathbb{S}^{d{-}1}}\!
  \!r^d K(r\xi)\dr\dH{\xi}{d-1} \int_{\T^d}|\nabla \mathbbm{1}_\Omega|=\int_{\R^d}|x|K(x)\dx \int_{\T^d}|\nabla \mathbbm{1}_\Omega|,
\end{align*}
where for the first inequality in the last line we reason as follows. Letting $v\in C_c^{\infty}(\T^d)$ and denoting $w(t)=v(x+r\sqrt{h}t\xi)$, we have
\begin{equation*}
    v(x+r\sqrt{h}\xi)-v(x)=\int_0^1 w'(t) \dt=r\sqrt{h}\int_0^1\nabla v(x+r\sqrt{h}t \xi){\cdot}  \xi \dt , 
\end{equation*}
and hence, by Fubini's theorem,
\begin{align}
    \frac{1}{r\sqrt{h}}\int_{\Omega}\left| v(x{+}\sqrt{h}\xi){-}v(x) \right|\dx  &\leq \int_0^1\int_{\Omega}|\nabla v(x{+}r\sqrt{h}t\xi)|\dx \dt\notag\\
    &= \int_0^1\int_{\Omega {+} r\sqrt{h}t\xi}|\nabla v(x)|\dx \dt\notag\\
    &\leq\int_0^1\int_{\T^d}|\nabla v(x)|\dx \dt  = \int_{\T^d}|\nabla v|\label{domination with perimeter}. 
\end{align}
Then the result for $\mathbbm{1}_\Omega=v\in BV(\T^d)$ follows by approximation of $v\in BV(\T^d)$ with functions $v_{\varepsilon}\in  C_c^{\infty}(\T^d) $ such that $v_{\varepsilon}\rightarrow v$ in $L^1(\T^d)$ and $\lim_{\varepsilon\rightarrow \infty} \int_{\T^d}|\nabla v_{\varepsilon}|\dx =\int_{\T^d}|\nabla v|$. 

Hence, combining \eqref{ineq term R^d minus A_x } with \eqref{ineq compactness 1}, we find
    \begin{align*}
 \int_{\Omega}|\nabla (J_h*u_h)|\dx  &\stackrel{\eqref{ineq compactness 1},\eqref{ineq term R^d minus A_x }}{\leq} 
 \frac{2 c}{\sqrt{h}} \int_{\Omega}(1-u_h) K_{4h} * u_h \dx +c \leq c(1+E_{4h}(u_h))\\
 &\stackrel{\eqref{eq: lemma gamma conv 2} }{\leq} c(1+(1+2\sqrt{h})E_h(u_h))\leq c,
    \end{align*}
    where we used Lemma \ref{lemma gamma conv 2}) with $N=2$.
    Thus, combining this with $|J_h*u_h|\leq 1$, by compactness in $BV$ we have (up to a non-relabeled subsequence) that $(J_h*u_h)_h$ converges in $L^1(\Omega)$. Moreover, we also have 
    \begin{align}
        \int_{\Omega}\left|J_h * u_h-u_h\right| \dx  &\stackrel{\eqref{ineq 2}}{\leq} \int_{\R^d} J_h(y) \int_{\Omega}|u_h(x+y)-u_h(x)| \dx  \dy\notag\\
        &\leq c_J\int_{\R^d} K_h(y) \int_{\Omega}|u_h(x+y)-u_h(x)| \dx  \dy\notag\\
        &\stackrel{\eqref{ineq 1},\eqref{ineq term R^d minus A_x }}{\leq}2c\int_{\Omega}(1-u_h) K_{h} * u_h \dx +c\sqrt{h} \leq c\sqrt{h}\label{compactness eq 1},
    \end{align} 
    so that also $(u_h)_h$ converges in $L^1(\Omega)$ up to a non-relabeled subsequence. 
Let now $u_h\rightarrow u$ in $L^1(\Omega)$. As seen above, we also have $J_h*u_h\rightarrow u$ in $L^1(\Omega)$ and $\nabla J_h*u_h$ bounded in $L^1(\Omega)$ uniformly in $h$, so that $u\in BV(\Omega)$. Finally, by
    \begin{align}
        \!\int_{\Omega}\! u_h(1{-}u_h) \dx  \stackrel{\eqref{ineq 3}}{\leq} \int_{\Omega}(1{-}u_h) J_h {*} u_h \dx {+}\!\int_{\Omega}\!\left|J_h {*} u_h{-}u_h\right| \dx \stackrel{   \eqref{eq compactness distance characteristic function preliminary 1},\eqref{compactness eq 1}}{\leq} c\sqrt{h},\label{eq compactness distance characteristic function}
    \end{align}
    we get that $u\in BV(\Omega;\{0,1\})$. Indeed, from the above inequality $\int_{\Omega} u_h^2$ converges, so $(u_h)_h$ converges both in $L^1$ and $L^2$ and thus the limit must be the same. Hence, $\int_{\Omega} u(1-u) \dx  =0$ and since both $u$ and $1-u$ are nonnegative, we have the result.
\end{proof}

We finally prove the pointwise convergence of the energy. 

\begin{proof}[Proof of Lemma \ref{lemma gamma conv 1}]
    We will show 
    \begin{align*}
        \lim_{h\rightarrow 0} \frac{1}{\sqrt{h}}\int_{\Omega}\gammatPV \,\tilde{v}K_h*v\dx=
            \frac{1}{2}\int_{ \R^d}\gammaPV \,\gamma(\nu)\left(|\nabla v|{+}|\nabla \tilde{v}|{-}|\nabla\mathbbm{1}_{\Omega}|\right) ,
    \end{align*}
    where $v,\tilde{v}\in BV(\Omega;[0,1])$, playing the roles of $u$ and $\mathbbm{1}_\Omega{-}u$ respectively,  are such that $v\tilde{v}=0$ almost everywhere. The corresponding limits for the $\gammatSP$- and $\gammatSV$-terms follow by analogous calculations.

    We first rewrite the left hand side in polar coordinates
\begin{align*}
\frac{1}{\sqrt{h}} \int_{\Omega}\gammatPV \ \tilde{v} K_h {* }v \dx  &=\frac{1}{\sqrt{h}} \int_{\mathbb{R}^d} K_h(y) \int_{\Omega} \gammatPV (x)\tilde{v}(x) v(x+y) \dx  \dy \\
& =\frac{1}{\sqrt{h}} \int_{\mathbb{R}^d} K(z) \int_{\Omega} \gammatPV (x)\tilde{v}(x) v(x+\sqrt{h} z) \dx  \dy 
\\
& =\int_0^{+\infty} r^d \frac{1}{r \sqrt{h}} \int_{\mathbb{S}^{d-1}} K(r \xi) \int_{\Omega} \gammatPV (x)\tilde{v}(x) v(x+r \sqrt{h} \xi) \dx  \dH{\xi}{d-1} \dr .
\end{align*}
In order to apply the dominated convergence theorem, we estimate
\begin{align}
    &\left| r^d \int_{\mathbb{S}^{d-1}} K(r \xi) \int_{\Omega} \frac{1}{r \sqrt{h}}\gammatPV (x)\tilde{v}(x) v(x+r \sqrt{h} \xi) \dx  \dH{\xi}{d-1}\right|\notag\\
    &\quad \quad\stackrel{\eqref{condition space gammat bounds}}{\leq} \widetilde{C}r^d \int_{\mathbb{S}^{d-1}} K(r \xi) \int_{\Omega} \frac{1}{r \sqrt{h}}\left|\tilde{v}(x) v(x+r \sqrt{h} \xi)\right| \dx  \dH{\xi}{d-1}\notag\\
    &\quad \quad=\widetilde{C}r^d \int_{\mathbb{S}^{d-1}} K(r \xi) \int_{\Omega} \frac{1}{r \sqrt{h}}\left|\tilde{v}(x) \left(v(x+r \sqrt{h} \xi)-v(x)\right)\right| \dx  \dH{\xi}{d-1}\notag\\
    &\quad \quad\leq \widetilde{C}r^d \int_{\mathbb{S}^{d-1}} K(r \xi) \int_{\Omega} \frac{1}{r \sqrt{h}}\left| v(x+r \sqrt{h} \xi)-v(x)\right| \dx  \dH{\xi}{d-1}  \notag \\ 
    &\quad \quad\stackrel{\eqref{domination with perimeter}}{\leq} \widetilde{C}r^d \int_{\mathbb{S}^{d-1}} K(r \xi) \dH{\xi}{d-1} \int_{\Omega}|\nabla v|,\label{dominating function 1}
\end{align}
where we used that $\tilde{v}v=0$ and $\tilde{v}\in[0,1]$ almost everywhere in the second equality and in the third inequality, respectively. The right hand side in \eqref{dominating function 1} is integrable, i.e., 
\begin{align*}
    \int_{\Omega}|\nabla v|\int_{0}^{\infty}\int_{\mathbb{S}^{d-1}} r^d K(r \xi) \dH{\xi}{d-1} \dr\leq c\int_{\R^d} |x| K(x)\dx <\infty,
\end{align*}
since $v\in BV(\Omega)$ and  $\int_{\R^d} |x| K(x)\dx \leq c$. 
For the pointwise limit we want to compute
\begin{equation*}
    \lim_{h\rightarrow 0}\int_{\mathbb{S}^{d-1}} r^d K(r \xi) \int_{\Omega} \frac{1}{r \sqrt{h}} \gammatPV(x) \tilde{v}(x) v(x+r \sqrt{h} \xi) \dx  \dH{\xi}{d-1}.
\end{equation*}
To this end, we apply again the dominated convergence theorem. 
To find a dominating function we write as above
\begin{equation*}
    \left|r^d K(r \xi) \int_{\Omega} \frac{1}{r \sqrt{h}}\gammatPV (x)\tilde{v}(x) v(x+r \sqrt{h} \xi) \dx  \right|\leq \widetilde{C} r^d K(r \xi) \int_{\Omega}|\nabla v|,
\end{equation*}
where $r^d K(r \xi) \in L^1(\mathbb{S}^{d-1})$. 
Then, for fixed $r>0$, let us compute the pointwise limit
\begin{align*}
    \lim_{h\rightarrow 0}\int_{\Omega} \frac{1}{r \sqrt{h}}\gammatPV (x)\tilde{v}(x) v(x+r \sqrt{h} \xi) \dx
\end{align*}
or, equivalently by the symmetry of $K$,
\begin{equation*}
    \lim_{h\rightarrow 0}\frac{1}{2}\int_{\Omega}\frac{1}{ \sqrt{h}}\gammatPV (x)\tilde{v}(x)\left(v(x+ \sqrt{h} \xi)+v(x- \sqrt{h} \xi) \right)\dx.
\end{equation*}
In order to prove the theorem, it is enough to show that  
\begin{align}
    &\lim_{h\rightarrow 0} \!\int_{\Omega}\!\frac{1}{ \sqrt{h}}\gammatPV (x)\tilde{v}(x)\left(v(x{+} \sqrt{h} \xi){+}v(x{-} \sqrt{h} \xi) \right) \dx\notag\\  
    &\quad=\frac{1}{2}\!\int_{\R^d}\!\gammatPV \left(|\xi{\cdot} \nabla v|{+}|\xi{\cdot} \nabla \tilde{v}|{-}|\xi{\cdot }\nabla (v{+}\tilde{v})|\right) =\frac{1}{2}\!\int_{\R^d}\!\gammatPV |\xi{\cdot}\nu|\left(| \nabla v|{+}| \nabla \tilde{v}|{-}|\nabla (v{+}\tilde{v})|\right), \label{pointwise convergence claim}
\end{align}
where $|\xi \cdot \nabla v|=|\xi \cdot \nu||\nabla v|$ with $\nu=-\frac{\nabla v}{|\nabla v|}$ denoting the measure-theoretic outer normal, which exists by Besicovitch differentiation theorem for measures.
 We claim that it is enough to show, for any $w,\tilde{w}\in BV(\R;\{0,1\})$ with $w\tilde{w}=0$ almost everywhere and $\widehat{\gamma}\in C_b(\R)$, 
\begin{equation}\label{pointwise convergence reduced 1d}
    \lim_{h\rightarrow 0 }\frac{1}{\sqrt{h}}\int_{\R}\widehat{\gamma}(t)\tilde{w}(t)\left(w(t{+}{\sqrt{h}}){+}w(t{-}{\sqrt{h}})\right)\dt=\frac{1}{2}\int_{\R}\widehat{\gamma}\left(\left|\frac{dw}{dt}\right|{+}\left|\frac{d\tilde{w}}{dt}\right|{-}\left|\frac{d(w{+}\tilde{w})}{dt}\right|\right).
\end{equation}
Indeed, first notice that, by symmetry of \eqref{pointwise convergence claim}, we can assume without loss of generality $\xi=e_d$. Let $w=v(x',\cdot)$, $\tilde{w}=\tilde{v}(x',\cdot)$, and $\widehat{\gamma}=\gammatPV (x',\cdot)$, where $x'=(x_1,\dots, x_{d-1})$. By \cite[Thm. 3.103]{AmbrosioFuscoPallara}, if $v\in BV(\R^d)$, then $w=v(x',\cdot)\in BV(\R)$ for almost every $x'\in \R^{d-1}$ and 
\begin{equation}\label{rel Dw and Dv}
    \int_{\R^{d-1}}\int_{\R}\widehat{\gamma}\left|\frac{dw}{dt}\right|\dx '=\int_{\R^d}\gammatPV |e_d{\cdot}\nabla v|\dx .
\end{equation}
Hence, \eqref{pointwise convergence claim} follows by \eqref{pointwise convergence reduced 1d}, \eqref{rel Dw and Dv}, and the dominated convergence theorem.

Let us thus show \eqref{pointwise convergence reduced 1d}.
Since $w,\tilde{w}\in BV(\R;\{0,1\})$ have a finite number of jumps between $0$ and $1$, which we denote by $t_1<\dots<t_n$ and $\tilde{t}_1<\dots<\tilde{t}_{\tilde{n}}$, respectively. Let us also denote as $\bar{t}_1<\dots<\bar{t}_N$ the ordered set of non repeated jump points of both $w$ and $\tilde{w}$, where $N=n+ \tilde{n} -\#(\{t_1,\dots,t_n\}\cap \{\tilde{t}_1,\dots,\tilde{t}_{\tilde{n}}\})$ and let
\begin{equation*}
    \tilde{I}\coloneqq \left\{i\in\{1,\dots,N-1\}\mid \tilde{w}\neq 0 \text{ in } (\bar{t}_i,\bar{t}_{i+1})\right\},
\end{equation*}
\begin{equation*}
    I\coloneqq \left\{i\in\{1,\dots,N\}\mid \bar{t}_i\in \{t_1,\dots,t_n\}\cap \{\tilde{t}_1,\dots,\tilde{t}_{\tilde{n}}\}\right\}.
\end{equation*}
Note that by definition we have $\tilde{w}\equiv 1$ in $\cup_{i\in \tilde{I}}(\bar{t}_{i},\bar{t}_{i+1})$ and $\tilde{w}\equiv 0$ in $\R\setminus \cup_{i\in \tilde{I}}(\bar{t}_{i},\bar{t}_{i+1})$.  
Hence,
\begin{align*}
    \frac{1}{\sqrt{h}}\!\!\int_{\R}\!\!\widehat{\gamma}(t)\tilde{w}(t)\left(w(t{+} \sqrt{h}){+}w(t{-} \sqrt{h} )\right){\rm d}t
    =\frac{1}{\sqrt{h}}\sum_{i\in \tilde{I}} \!\int_{\bar{t}_i}^{\bar{t}_{i+1}}\!\!\!\!\widehat{\gamma}(t)\left(w(t{+} \sqrt{h}){+}w(t{-} \sqrt{h} )\right){\rm d}t,
\end{align*}
where we used that $\tilde{w}$ is compactly supported. Moreover, we can also assume $\sqrt{h}<\min_{i\neq j}|\bar{t}_i-\bar{t}_j |$. Thus, since $\tilde{w}w=0$ almost everywhere, we deduce that for $t\in (\bar{t}_{i},\bar{t}_{i+1})$, $i\in\tilde{I} $,   $w(t+\sqrt{h})\neq 0$ only if $t\in (\bar{t}_{i+1}-\sqrt{h},\bar{t}_{i+1})$ and similarly $w(t-\sqrt{h})\neq 0$ only if $t\in (\bar{t}_{i},\bar{t}_{i}+\sqrt{h})$.
Consequently, we have that
\begin{align*}
    \frac{1}{\sqrt{h}}\int_{\R}\!\!\widehat{\gamma}(t)\tilde{w}(t)\left(w(t{+} \sqrt{h}){+}w(t{-} \sqrt{h} )\right){\rm d}t
    &=\frac{1}{\sqrt{h}}\sum_{i\in \tilde{I}} \!\int_{\bar{t}_{i+1}-\sqrt{h'}}^{\bar{t}_{i+1}}\!\widehat{\gamma}(t)w(t{+} \sqrt{h})\dt\\
    &\quad+\frac{1}{\sqrt{h}}\sum_{i\in \tilde{I}} \!\int_{\bar{t}_{i}}^{\bar{t}_{i}+\sqrt{h}}\!\widehat{\gamma}(t)w(t{-} \sqrt{h})\dt.
\end{align*}
Now notice that $w(t+\sqrt{h})$ and $w(t-\sqrt{h})$ are non-zero if and only if $t_i$ or $t_{i+1}$ are joint jumps of $w$ and $\tilde{w}$. Indeed, looking without loss of generality at the case of $w(t+\sqrt{h})$, we see that,
$w\neq 0$ in $(\bar{t}_{i+1},\bar{t}_{i+1}+\sqrt{h})$ implies $\tilde{w}=0$ in $(\bar{t}_{i+1},\bar{t}_{i+1}+\sqrt{h})$, i.e., $\bar{t}_{i+1}$ is a jump point for both $w$ and $\tilde{w}$.
Hence, by Lebesgue differentiation theorem and by the definition of $I$, we get 
\begin{equation*}
    \lim_{h\rightarrow 0}\frac{1}{\sqrt{h}}\int_{\R}\!\!\widehat{\gamma}(t)\tilde{w}(t)\left(w(t{+} \sqrt{h}){+}w(t{-} \sqrt{h} )\right){\rm d}t=\sum_{i\in I}\widehat{\gamma}(\bar{t}_{i}).
\end{equation*}
On the other hand, we have,
\begin{equation*}
   \int_{\R}\widehat{\gamma}(t)\left|\frac{dw}{dt}\right|=\sum_{i=1}^{n}\widehat{\gamma}(t_i),\quad \int_{\R}\widehat{\gamma}(t)\left|\frac{d\tilde{w}}{dt}\right|=\sum_{i=1}^{\tilde{n}}\widehat{\gamma}(\tilde{t}_i),
\end{equation*}
\begin{equation*}
  \int_{\R}\widehat{\gamma}(t)\left|\frac{d(w+\tilde{w})}{dt}\right|=\sum_{i=1}^{n}\widehat{\gamma}(t_i)+\sum_{i=1}^{\tilde{n}}\widehat{\gamma}(\tilde{t}_i)-2\sum_{i\in I}\widehat{\gamma}(\bar{t}_i),
\end{equation*}
and we have thus obtained \eqref{pointwise convergence claim}.
We have shown that 
\begin{align*}
     \lim_{h\rightarrow 0} \frac{1}{\sqrt{h}}&\int_{\Omega}\gammatPV \tilde{v}K_h*v\dx 
     \\&=\frac{1}{4}\int_{ \R^d}\left(\!\int_0^{+\infty}\!\!\! r^d  \!\!\int_{\mathbb{S}^{d-1}}\!\! |\xi\cdot\nu| K(r \xi) \dH{\xi}{d-1} \dr\right)\gammatPV  \left(| \nabla v|{+}| \nabla \tilde{v}|{-}|\nabla (v{+}\tilde{v})|\right) \\
     &\stackrel{\eqref{gamma from kernel}}{=}\frac{1}{2}\int_{ \R^d}\gammatPV \gamma(\nu)\left(| \nabla v|{+}| \nabla \tilde{v}|{-}|\nabla \mathbbm{1}_{\Omega}|\right) ,
\end{align*}
which concludes the proof.
\end{proof}

\section*{Acknowledgments}
AC is supported by the Austrian Science Fund (FWF) projects 10.55776/F65,  10.55776/I5149,
10.55776/P32788,     and   10.55776/I4354,
as well as by the OeAD-WTZ project CZ 09/2023. 
Part of this research was funded by the Mobility Fellowship of the University of Vienna and conducted during a visit to Kyoto University, whose warm hospitality is gratefully acknowledged.
The work of KS was supported by JSPS KAKENHI Grant Number 19K03634.
%For open-access purposes, the authors have applied a CC BY public copyright license to any author-accepted manuscript version arising from this submission. 

\end{document}